%% file: header.tex
\input{arxiv}

\expandafter\ifx\csname c@nestcount\endcsname\relax\newcounter{nestcount}\fi
\global\addtocounter{nestcount}{1}
\ifnum\value{nestcount}>1\else

\documentclass[12pt,reqno]{amsart}

\usepackage{a4wide}
\usepackage{amsmath}
\usepackage{amssymb}
\usepackage{amsxtra}
\usepackage{amsthm}
\usepackage{mathtools}
\usepackage{extpfeil}
\usepackage{enumitem}
\usepackage[mathscr]{eucal}
\usepackage{graphicx}
\usepackage{asymptote}
\usepackage{proof}
\usepackage{ifnextok}
\usepackage[T1]{fontenc}
\usepackage{xstring}

\usepackage[all]{xy}
\SelectTips{cm}{}
\SilentMatrices
\ifx\pdfoutput\undefined
  \xyoption{dvips}
\fi

\usepackage[silent]{fixme}
\usepackage[unicode,final]{hyperref}
\usepackage{slashed}

\usepackage{mathrsfs}
\usepackage{xr}

\ifarxiv
  \fxsetup{
    status=final
  }
\else
  \usepackage[inline]{showlabels}
  \fxsetup{
      status=draft,
      layout={marginclue,noinline}
  \fxsetface{inline}{\itshape}
\fi

  \ifarxiv
    \input{macros}
  \else
    \input{../common/macros}
  \fi

\setlist{}
\setenumerate{leftmargin=*,labelindent=\parindent}
\setitemize{leftmargin=*,labelindent=0.5\parindent}
\setdescription{leftmargin=1em}

\swapnumbers

\theoremstyle{plain}

\theoremstyle{definition}

\theoremstyle{remark}

\makeatletter
\let\c@equation\c@thm
\makeatother
\numberwithin{equation}{subsection}

\raggedbottom

\ifarxiv
  \begin{asydef}
    import "../common/adjunction" as adjunction;
  \end{asydef}
\else
  \begin{asydef}
    import "adjunction" as adjunction;
  \end{asydef}
\fi

\ifarxiv
  \input{foundations_lab}
  \input{cohadj_lab}
  \input{reedy_lab}
  \input{complete_lab}
  \input{equipment_lab}
  \input{yoneda_lab}
\else
  \externaldocument[found:]{../foundations/all}
  \externaldocument[cohadj:]{../cohadj/all}
  \externaldocument[reedy:]{../reedy/reedy}
  \externaldocument[complete:]{../completeness/completeness}
  \externaldocument[equipment:]{../equipment/all}
  \externaldocument[yoneda:]{../yoneda/all}
\fi

\newcommand{\extRef}[3]{%
  {\protect\IfBeginWith{#3}{itm:}{}{#2.}}\ref*{#1:#3}}

\newcommand{\refI}{\extRef{found}{I}}

\newcommand{\refIV}{\extRef{yoneda}{IV}}

\title{Working Title}

\author[Riehl]{Emily Riehl}
\address{
  Department of Mathematics \\
Johns Hopkins University \\
Baltimore, MD 21218\\
  USA
}
\email{eriehl@math.jhu.edu}

\author[Verity]{Dominic Verity}
\address{
  Centre of Australian Category Theory \\
  Macquarie University \\
  NSW 2109 \\
  Australia
}
\email{dominic.verity@mq.edu.au}

\date{\today}


\setcounter{tocdepth}{1}
\expandafter\ifx\csname preamble\endcsname\relax\else\preamble\fi

\begin{document}

  \ifpdf
  \DeclareGraphicsExtensions{.pdf, .jpg, .tif}
  \else
  \DeclareGraphicsExtensions{.eps, .jpg}
  \fi

  \expandafter\ifx\csname abstracttext\endcsname\relax\else
  \begin{abstract}
    \abstracttext
  \end{abstract}
  \fi

  \maketitle

\tableofcontents

\fi

%% file: arxiv.tex
\newif\ifarxiv

\arxivtrue

%% file: macros.tex

\input{arxiv}

\newcommand{\ifundef}[1]{\expandafter\ifx\csname#1\endcsname\relax}



\DeclareMathAlphabet{\mathbbe}{U}{bbold}{m}{n}

\makeatletter

\def\re@DeclareMathSymbol#1#2#3#4{%
    \let#1=\undefined
    \DeclareMathSymbol{#1}{#2}{#3}{#4}}

\ifundef{Top}
  \DeclareSymbolFont{tcSyC}{U}{txsyc}{m}{n}
  \SetSymbolFont{tcSyC}{bold}{U}{txsyc}{bx}{n}
  \DeclareFontSubstitution{U}{txsyc}{m}{n}

  \re@DeclareMathSymbol{\Top}{\mathord}{tcSyC}{120}
  \re@DeclareMathSymbol{\Bot}{\mathord}{tcSyC}{121}
\fi

\ifundef{righthalfcup}
  \DeclareFontFamily{U}{MnSymbolC}{}
  \DeclareSymbolFont{mnSyC}{U}{MnSymbolC}{m}{n}
  \SetSymbolFont{mnSyC}{bold}{U}{MnSymbolC}{b}{n}
  \DeclareFontShape{U}{MnSymbolC}{m}{n}{
      <-6>  MnSymbolC5
     <6-7>  MnSymbolC6
     <7-8>  MnSymbolC7
     <8-9>  MnSymbolC8
     <9-10> MnSymbolC9
    <10-12> MnSymbolC10
    <12->   MnSymbolC12}{}
  \DeclareFontShape{U}{MnSymbolC}{b}{n}{
      <-6>  MnSymbolC-Bold5
     <6-7>  MnSymbolC-Bold6
     <7-8>  MnSymbolC-Bold7
     <8-9>  MnSymbolC-Bold8
     <9-10> MnSymbolC-Bold9
    <10-12> MnSymbolC-Bold10
    <12->   MnSymbolC-Bold12}{}
  
  \re@DeclareMathSymbol{\righthalfcup}{\mathord}{mnSyC}{184}
  \re@DeclareMathSymbol{\lefthalfcap}{\mathord}{mnSyC}{185}
\fi

\DeclareFontFamily{U}{MnSymbolA}{}
\DeclareSymbolFont{mnSyA}{U}{MnSymbolA}{m}{n}
\SetSymbolFont{mnSyA}{bold}{U}{MnSymbolA}{b}{n}
\DeclareFontShape{U}{MnSymbolA}{m}{n}{
    <-6>  MnSymbolA5
   <6-7>  MnSymbolA6
   <7-8>  MnSymbolA7
   <8-9>  MnSymbolA8
   <9-10> MnSymbolA9
  <10-12> MnSymbolA10
  <12->   MnSymbolA12}{}
\DeclareFontShape{U}{MnSymbolA}{b}{n}{
    <-6>  MnSymbolA-Bold5
   <6-7>  MnSymbolA-Bold6
   <7-8>  MnSymbolA-Bold7
   <8-9>  MnSymbolA-Bold8
   <9-10> MnSymbolA-Bold9
  <10-12> MnSymbolA-Bold10
  <12->   MnSymbolA-Bold12}{}

\re@DeclareMathSymbol{\twoheadedswarrow}{\mathord}{mnSyA}{30}

\makeatother




\newcommand{\mlaux}[3]{\setbox0=\hbox{$\mathsurround=0pt #2{#3}$}%
  \dimen0=\dp0\advance\dimen0 by \ht0\lower#1\dimen0\box0}

\newcommand{\makellapm}[2]{\hbox to 0pt{\hss$\mathsurround=0pt #1{#2}$}}

\newcommand{\makerlapm}[2]{\hbox to 0pt{$\mathsurround=0pt #1{#2}$\hss}}

\newcommand{\makelapm}[2]{\hbox to 0pt{\hss$\mathsurround=0pt #1{#2}$\hss}}

\newcommand{\makeushort}[3]{%
	\setbox0=\hbox{$\mathsurround=0pt #2{#3}$}%
	\hbox to 1\wd0{\hss\underbar{\hbox to #1\wd0{\hss\box0\hss}}\hss}}

\def\makebigger#1#2#3{\scalebox{#1}{$\mathsurround=0pt #2{#3}$}}
\def\bigger#1#2{{\relax\mathpalette{\makebigger{#1}}{#2}}}

\def\scaleuphalf{1.0954}
\def\scaleupone{1.2}





\makeatletter

\def\newmop{\@ifstar{\@newmop m}{\@newmop o}}
\def\@newmop#1{\@ifnextchar[{\@@newmop #1}{\@@@newmop #1}}
\def\@@newmop#1[#2]{\@declmathop #1#2}
\def\@@@newmop#1#2{\expandafter\@declmathop\expandafter #1\csname #2\endcsname{#2}}

\makeatother

\newdir{ |}{{}*!/-5pt/@{|}}

\newmop{im}
\newmop{coim}
\newmop{dom}
\newmop{cod}
\newmop{id}
\newmop{Map}

\newmop{obj}
\newmop{arr}
\newmop{sq}
\newmop{norm}

\newmop{el}

\newmop{ev}


\newmop{Ext}
\newmop{icon}
\newmop{pbk}
\newmop{icom}

\newmop{cell}
\newmop{cof}
\newmop{fib}

\newmop{diag}

\newmop*{colim}
\newmop*{holim}

\newmop[\lan]{lan}
\newmop[\ran]{ran}


\makeatletter
\newcommand{\rotatemath}[2]{\rotatebox[origin=c]{180}{$\m@th #1{#2}$}}
\makeatother

\newmop[\pr]{pr}
\newmop[\dm]{d} 

\newmop[\cpr]{in}


\newcommand{\pocorner}{\hbox to 8pt{{\vrule height8pt depth0pt width0.5pt}%
    \vbox to 8pt{{\hrule height0.5pt width7.5pt depth0pt}\vfill}}}

\newcommand{\pbcorner}{\vbox to 0pt{\kern 4pt\hbox to 0pt{\kern 4pt%
      \vbox{{\hrule height0.5pt width7.5pt depth0pt}}%
      {\vrule height8pt depth0pt width0.5pt}\hss}\vss}}


\newmop{cls}



\newcommand{\category}[1]{\underline{\smash[b]{\text{\rm{#1}}}}}


\makeatletter

\def\Del@Sym{{\bigger\scaleuphalf{\mathbbe{\Delta}}}}

\def\del@fn{\futurelet\del@next}
\def\del@dn{\def\del@next}

\def\parsedel@{%
  \ifx +\del@next \del@dn+{\Del@Sym_{\mathord{+}}}%
  \else \del@dn {\del@fn\parsedel@@}%
  \fi\del@next}

\def\parsedel@@{%
  \ifx\space@\del@next \expandafter\del@dn\space{\del@fn\parsedel@@}%
  \else\ifx [\del@next \del@dn[{\del@fn\parsedel@@@}%
  \else\ifx _\del@next \del@dn{\Delta}%
  \else\ifx ^\del@next \del@dn{\Delta}%
  \else \del@dn{\Del@Sym}%
  \fi\fi\fi\fi\del@next}

\def\parsedel@@@{%
  \ifx\space@\del@next \expandafter\del@dn\space{\del@fn\parsedel@@@}%
  \else\ifx t\del@next \del@dn t{\Del@Sym_\infty\del@fn\parsedel@@@@}%
  \else\ifx b\del@next \del@dn b{\Del@Sym_{-\infty}\del@fn\parsedel@@@@}%
  \else \del@dn{\errmessage{unexpected modifier}}%
  \fi\fi\fi\del@next}

\def\parsedel@@@@{%
  \ifx\space@\del@next \expandafter\del@dn\space{\del@fn\parsedel@@@@}%
  \else\ifx ]\del@next \del@dn]{}%
  \else \del@dn{\errmessage{expecting close of option block}}%
  \fi\fi\del@next}

\def\Del{\del@fn\parsedel@}

\makeatother



\newcommand{\qCat}{\category{qCat}}
\newcommand{\Adj}{\category{Adj}}
\newcommand{\Mnd}{\category{Mnd}}



\newcommand{\dmod}[3]{\xymatrix@=1.25em{{#2} \ar[r]|\mid^{ {#1}} & {#3}}}


\newcommand{\pbshape}{{\mathord{\bigger\scaleupone\righthalfcup}}}






\newmop{dec}
\newmop{fatdec}
\newmop{slc}
\newmop{fatslc}

\newmop{ir} 
\newmop{incl}





\def\reedyfilt#1_#2{#1_{\leq #2}}
\newmop{sk}
\newmop{cosk}
\newmop{res}



\newmop{map}
\newmop{Ho}

\newmop[\pth]{path}
\newmop{cyl}



\newmop[\const]{c}
  



\newmop{coll}
\newmop{wgt}




\newdir{ >}{{}*!/-7pt/@{>}}
\newdir{u(}{{}*!/-4pt/@^{(}}
\newdir{d(}{{}*!/-4pt/@_{(}}
\newdir{|>}{%
  !/4.5pt/@{|}*:(1,-.2)@^{>}*:(1,+.2)@_{>}*+@{}}


\makeatletter
\def\makeslashed#1#2#3#4#5{#1{\mathpalette{\sla@{#2}{#3}{#4}}{#5}}}

\def\@mathlower#1#2#3{\setbox0=\hbox{$\m@th#2#3$}\lower#1\ht0\box0}
\def\mathlower#1#2{\mathpalette{\@mathlower{#1}}{#2}}
\makeatother




\makeatletter

\def\tens@fn{\futurelet\tens@next}
\def\tens@dn{\def\tens@nextcont}
\newtoks\tens@toks
\def\addtotens@toks#1{\tens@toks=\expandafter{\the\tens@toks#1}}

\def\parsetens@@{%
    \ifx\space@\tens@next \expandafter\tens@dn\space{\tens@fn\parsetens@@}%
    \else\ifx ^\tens@next \tens@dn ^##1{\parsetens@procsep^\addtotens@toks{##1}%
      \tens@fn\parsetens@@}%
    \else\ifx _\tens@next \tens@dn _##1{\parsetens@procsep_\addtotens@toks{##1}%
      \tens@fn\parsetens@@}%
    \else\tens@dn{\ifx *\tens@last \else\addtotens@toks\egroup\fi\the\tens@toks}%
    \fi\fi\fi\tens@nextcont}

\def\parsetens@procsep#1{%
  \ifx *\tens@last \addtotens@toks{#1}\addtotens@toks\bgroup%
  \else\ifx \tens@last\tens@next \addtotens@toks,%
  \else \addtotens@toks\egroup\addtotens@toks\bgroup%
    \addtotens@toks\egroup\addtotens@toks{#1}\addtotens@toks\bgroup%
  \fi\fi\let\tens@last\tens@next}

\newcommand{\tn}[1]{\let\tens@last=*\tens@toks={#1}\tens@fn\parsetens@@}

\makeatother


\def\adjdisplay#1-|#2:#3->#4.{{%
    \xymatrix@R=0em@!C=2.5em{%
      *+[l]{#3} \ar@/_0.55pc/[rr]_-{#2} & {\bot} &
      *+[r]{#4}\ar@/_0.55pc/[ll]_-{#1}}}}

\def\adjdisplaytwo#1-|#2:#3->#4.{{%
\xymatrix@=1.2em{
      {#3}\ar@/_1.5ex/[rr]_-{#2}^-{}="one"
      & & {#4}
      \ar@/_1.5ex/[ll]_-{#1}^-{}="two" 
      \ar@{}"one";"two"|{\bot}
    }}}

\def\tripleadjdisplay#1-|#2-|#3:#4->#5.{{%
\xymatrix@=2.4em{ 
{#4}\ar[r]|{#2} &
{#5} \ar@/_3ex/[l]_{#1}^{\bot} \ar@/^3ex/[l]_{\bot}^{#3}}
}}

\def\adjinline#1-|#2:#3->#4.{{#1}\dashv{#2}:#3\to #4}

\newcommand{\pent}[1]{
  \xybox{
    \POS (0,-15)*+{\a}="0", 
         (-14,-5)*+{\b}="1", 
         (-9,12)*+{\c}="2", 
         (9,12)*+{\d}="3", 
         (14,-5)*+{\e}="4"
    \POS"0" \ar "1"^{\labelstyle \ab}|{}="01"
    \POS"1" \ar "2"^{\labelstyle \bc}|{}="12"
    \POS"2" \ar "3"^{\labelstyle \cd}|{}="23"
    \POS"3" \ar "4"^{\labelstyle \de}|{}="34"
    \POS"0" \ar "4"_{\labelstyle \ae}|{}="04"
    \ifcase #1
    \POS"0" \ar "2"|{\labelstyle \ac}="02"
    \POS"0" \ar "3"|{\labelstyle \ad}="03"
    \POS"02";"1"**{}, ?(0.3) \ar@{=>} ?(0.7)^{\labelstyle \abc}
    \POS"03";"2"**{}, ?(0.25) \ar@{=>} ?(0.5)_{\labelstyle \acd}
    \POS"04";"3"**{}, ?(0.2) \ar@{=>} ?(0.4)_{\labelstyle \ade}
    \or
    \POS"1" \ar "3"|{\labelstyle \bd}="13"
    \POS"1" \ar "4"|{\labelstyle \be}="14"
    \POS"13";"2"**{}, ?(0.3) \ar@{=>} ?(0.7)_{\labelstyle \bcd}
    \POS"14";"3"**{}, ?(0.25) \ar@{=>} ?(0.5)_{\labelstyle \bde}
    \POS"04";"1"**{}, ?(0.25) \ar@{=>} ?(0.5)_{\labelstyle \abe}
    \or
    \POS"2" \ar "4"|{\labelstyle \ce}="24"
    \POS"0" \ar "2"|{\labelstyle \ac}="02"
    \POS"02";"1"**{}, ?(0.3) \ar@{=>} ?(0.7)^{\labelstyle \abc}
    \POS"04";"2"**{}, ?(0.2) \ar@{=>} ?(0.35)_{\labelstyle \ace}
    \POS"24";"3"**{}, ?(0.2) \ar@{=>} ?(0.6)^{\labelstyle \cde}
    \or
    \POS"1" \ar "3"|{\labelstyle \bd}="13"
    \POS"0" \ar "3"|{\labelstyle \ad}="03"
    \POS"04";"3"**{}, ?(0.2) \ar@{=>} ?(0.4)_{\labelstyle \ade}
    \POS"13";"2"**{}, ?(0.3) \ar@{=>} ?(0.7)_{\labelstyle \bcd}
    \POS"03";"1"**{}, ?(0.25) \ar@{=>} ?(0.5)^{\labelstyle \abd}
    \or
    \POS"2" \ar "4"|{\labelstyle \ce}="24"
    \POS"1" \ar "4"|{\labelstyle \be}="14"
    \POS"24";"3"**{}, ?(0.2) \ar@{=>} ?(0.6)^{\labelstyle \cde}
    \POS"04";"1"**{}, ?(0.25) \ar@{=>} ?(0.5)_{\labelstyle \abe}
    \POS"14";"2"**{}, ?(0.25) \ar@{=>} ?(0.5)^{\labelstyle \bce}
    \else\fi
  }
}

\newcommand{\pentofpent}[1]{
  \def\baselen{#1}
  \begin{xy}
    0;<\baselen,0mm>:
    *{\xybox{
        \POS(0,-4)*[o]{\pent 0}="zero"
        \POS(16,40)*[o]{\pent 3}="three"
        \POS(72,40)*[o]{\pent 1}="one"
        \POS(88,-4)*[o]{\pent 4}="four"
        \POS(44,-36)*[o]{\pent 2}="two"
        \ar@<1ex>"zero";"three"^-{\objectstyle\abcd}
        \ar@<1ex>"three";"one"^-{\objectstyle\abde}
        \ar@<1ex>"one";"four"^-{\objectstyle\bcde}
        \ar@<-1ex>"zero";"two"_-{\objectstyle\acde}
        \ar@<-1ex>"two";"four"_-{\objectstyle\abce}
        \ar@{=>}(44,-5);(44,+15)^{\objectstyle\abcde}
     }}
  \end{xy}
}


%% file: foundations_lab.tex
\newlabel{found:sec:background}{{2}{6}{Background on quasi-categories}{section.2}{}}
\newlabel{found:obs:size-conventions}{{2.0.1}{6}{size}{thm.2.0.1}{}}
\newlabel{found:subsect:simplicial.notation}{{2.1}{6}{Some standard simplicial notation}{subsection.2.1}{}}
\newlabel{found:ntn:simp.op}{{2.1.1}{6}{simplicial operators}{thm.2.1.1}{}}
\newlabel{found:ntn:simplicial-sets}{{2.1.2}{6}{(augmented) simplicial sets}{thm.2.1.2}{}}
\newlabel{found:rec:augmentation}{{2.1.3}{7}{augmentation}{thm.2.1.3}{}}
\newlabel{found:ntn:faces-by-vertices}{{2.1.5}{8}{faces of \protect {$\Del ^n$}}{thm.2.1.5}{}}
\newlabel{found:ntn:simplicial-hom-space}{{2.1.6}{8}{internal hom}{thm.2.1.6}{}}
\newlabel{found:rec:hty-category}{{2.2.2}{8}{the homotopy category}{thm.2.2.2}{}}
\newlabel{found:eq:homotopy-of-1-simplices}{{2.2.3}{9}{the homotopy category}{equation.2.2.3}{}}
\newlabel{found:rec:qmc-quasicat}{{2.2.5}{9}{the model category of quasi-categories}{thm.2.2.5}{}}
\newlabel{found:rec:leibniz}{{2.2.6}{9}{Leibniz constructions}{thm.2.2.6}{}}
\newlabel{found:rec:cart-modcat}{{2.2.8}{10}{cartesian model categories}{thm.2.2.8}{}}
\newlabel{found:obs:isofibration-closure}{{2.2.9}{10}{closure properties of isofibrations}{thm.2.2.9}{}}
\newlabel{found:defn:equivalences}{{2.3.1}{10}{isomorphisms in quasi-categories}{thm.2.3.1}{}}
\newlabel{found:rec:marked-prod-exp}{{2.3.4}{11}{products and exponentiation}{thm.2.3.4}{}}
\newlabel{found:rmk:equiv-markings}{{2.3.5}{12}{isomorphisms and markings}{thm.2.3.5}{}}
\newlabel{found:prop:joyal-special-horn}{{2.3.7}{12}{Joyal}{thm.2.3.7}{}}
\newlabel{found:rec:qmc-quasi-marked}{{2.3.8}{12}{the model structure of naturally marked quasi-categories}{thm.2.3.8}{}}
\newlabel{found:obs:nat-mark-homs}{{2.3.9}{13}{natural markings, internal homs, and products}{thm.2.3.9}{}}
\newlabel{found:lem:pointwise-equiv}{{2.3.10}{13}{pointwise isomorphisms are isomorphisms}{thm.2.3.10}{}}
\newlabel{found:eq:pointwise-equivalence-square}{{2.3.11}{13}{Isomorphisms and marked simplicial sets}{equation.2.3.11}{}}
\newlabel{found:cor:pointwise-equiv}{{2.3.12}{14}{}{thm.2.3.12}{}}
\newlabel{found:obs:marked-arrow-subcat}{{2.3.13}{14}{}{thm.2.3.13}{}}
\newlabel{found:subsec:join}{{2.4}{14}{Join and slice}{subsection.2.4}{}}
\newlabel{found:defn:join-dec}{{2.4.1}{15}{joins and d{\'e}calage}{thm.2.4.1}{}}
\newlabel{found:defn:slices}{{2.4.2}{15}{d{\'e}calage and slices}{thm.2.4.2}{}}
\newlabel{found:obs:slice-and-qcats}{{2.4.3}{16}{slices of quasi-categories}{thm.2.4.3}{}}
\newlabel{found:def:fat-join}{{2.4.4}{16}{fat join}{thm.2.4.4}{}}
\newlabel{found:eq:fat-join-def}{{2.4.5}{16}{fat join}{equation.2.4.5}{}}
\newlabel{found:eq:fat-join-cong}{{2.4.6}{17}{fat join}{equation.2.4.6}{}}
\newlabel{found:defn:fat-slices}{{2.4.7}{17}{fat slice}{thm.2.4.7}{}}
\newlabel{found:eq:join-comp-def}{{2.4.9}{17}{comparing join constructions}{equation.2.4.9}{}}
\newlabel{found:eq:tnm-def}{{2.4.10}{17}{comparing join constructions}{equation.2.4.10}{}}
\newlabel{found:prop:join-fatjoin-equiv}{{2.4.11}{18}{}{thm.2.4.11}{}}
\newlabel{found:lem:slices-quillen}{{2.4.12}{18}{}{thm.2.4.12}{}}
\newlabel{found:prop:slice-fatslice-equiv}{{2.4.13}{18}{slices and fat slices of a quasi-category are equivalent}{thm.2.4.13}{}}
\newlabel{found:rmk:map-slices}{{2.4.14}{19}{}{thm.2.4.14}{}}
\newlabel{found:sec:twocat}{{3}{19}{The 2-category of quasi-categories}{section.3}{}}
\newlabel{found:obs:simp-to-2-cats}{{3.1.2}{20}{}{thm.3.1.2}{}}
\newlabel{found:def:qCat-2}{{3.2.1}{21}{the 2-category of quasi-categories}{thm.3.2.1}{}}
\newlabel{found:eq:qCat2homdefn}{{3.2.2}{21}{the 2-category of quasi-categories}{equation.3.2.2}{}}
\newlabel{found:obs:pointwise-iso-reprise}{{3.2.3}{22}{pointwise isomorphisms are isomorphisms (reprise)}{thm.3.2.3}{}}
\newlabel{found:prop:qcat2closed}{{3.2.4}{22}{}{thm.3.2.4}{}}
\newlabel{found:defn:2-cat.of.all.simpsets}{{3.2.5}{23}{the 2-category of all simplicial sets}{thm.3.2.5}{}}
\newlabel{found:rmk:exp2functor}{{3.2.6}{23}{}{thm.3.2.6}{}}
\newlabel{found:lem:2-cat.equivs}{{3.2.8}{23}{}{thm.3.2.8}{}}
\newlabel{found:prop:equivsareequivs}{{3.2.9}{24}{}{thm.3.2.9}{}}
\newlabel{found:prop:equivsareequivs2}{{3.2.10}{24}{}{thm.3.2.10}{}}
\newlabel{found:subsec:weak-2-limits}{{3.3}{24}{Weak 2-limits}{subsection.3.3}{}}
\newlabel{found:defn:smothering}{{3.3.1}{24}{smothering functors}{thm.3.3.1}{}}
\newlabel{found:lem:smothering}{{3.3.2}{25}{fibres of smothering functors}{thm.3.3.2}{}}
\newlabel{found:defn:weak2limit}{{3.3.3}{25}{weak 2-limits in a 2-category}{thm.3.3.3}{}}
\newlabel{found:eq:cone-induced}{{3.3.4}{25}{weak 2-limits in a 2-category}{equation.3.3.4}{}}
\newlabel{found:lem:unique-weak-2-limits}{{3.3.5}{25}{}{thm.3.3.5}{}}
\newlabel{found:lem:weak-simplification}{{3.3.6}{26}{}{thm.3.3.6}{}}
\newlabel{found:eq:qCat.wl.comp.1}{{3.3.7}{26}{Weak 2-limits}{equation.3.3.7}{}}
\newlabel{found:prop:weak-cotensors}{{3.3.9}{27}{}{thm.3.3.9}{}}
\newlabel{found:eq:arrowinhA}{{3.3.10}{27}{Weak 2-limits}{equation.3.3.10}{}}
\newlabel{found:eq:liftedarrowinA2}{{3.3.11}{28}{Weak 2-limits}{equation.3.3.11}{}}
\newlabel{found:prop:weak-iso-cotensor}{{3.3.13}{28}{}{thm.3.3.13}{}}
\newlabel{found:prop:weak-homotopy-pullbacks}{{3.3.14}{29}{}{thm.3.3.14}{}}
\newlabel{found:def:comma-obj}{{3.3.15}{30}{comma objects}{thm.3.3.15}{}}
\newlabel{found:lem:comma-obj}{{3.3.16}{30}{}{thm.3.3.16}{}}
\newlabel{found:lem:comma-obj-maps}{{3.3.17}{30}{maps induced between comma objects}{thm.3.3.17}{}}
\newlabel{found:prop:weakcomma}{{3.3.18}{31}{}{thm.3.3.18}{}}
\newlabel{found:eq:commacat-comp}{{3.3.19}{31}{Weak 2-limits}{equation.3.3.19}{}}
\newlabel{found:obs:unpacking-weak-comma-objects}{{3.3.20}{31}{unpacking the universal property of weak comma objects}{thm.3.3.20}{}}
\newlabel{found:eq:weak-comma-prop}{{3.3.21}{31}{unpacking the universal property of weak comma objects}{equation.3.3.21}{}}
\newlabel{found:eq:standard-comma-pic}{{3.3.22}{32}{unpacking the universal property of weak comma objects}{equation.3.3.22}{}}
\newlabel{found:eq:comma-cone}{{3.3.23}{32}{unpacking the universal property of weak comma objects}{equation.3.3.23}{}}
\newlabel{found:eq:comma-ind-1cell-prop}{{3.3.24}{32}{unpacking the universal property of weak comma objects}{equation.3.3.24}{}}
\newlabel{found:eq:comma-ind-2cell-data}{{3.3.25}{32}{unpacking the universal property of weak comma objects}{equation.3.3.25}{}}
\newlabel{found:eq:comma-ind-2cell-compat}{{3.3.26}{33}{unpacking the universal property of weak comma objects}{equation.3.3.26}{}}
\newlabel{found:lem:1cell-ind-uniqueness}{{3.3.27}{33}{1-cell induction is unique up to isomorphism}{thm.3.3.27}{}}
\newlabel{found:subsec:slice-2cats-of-qcats}{{3.4}{33}{Slices of the category of quasi-categories}{subsection.3.4}{}}
\newlabel{found:defn:enriched-slice}{{3.4.1}{33}{enriching the slices of $\qCat $}{thm.3.4.1}{}}
\newlabel{found:eq:slice-hom-objects}{{3.4.2}{34}{enriching the slices of $\qCat $}{equation.3.4.2}{}}
\newlabel{found:obs:fibred-pushforward}{{3.4.3}{34}{pushforward}{thm.3.4.3}{}}
\newlabel{found:obs:fibred-pullback}{{3.4.4}{34}{pullback}{thm.3.4.4}{}}
\newlabel{found:defn:smothering-2-functor}{{3.4.6}{35}{smothering 2-functor}{thm.3.4.6}{}}
\newlabel{found:prop:slice-smothering-2-functor}{{3.4.7}{35}{}{thm.3.4.7}{}}
\newlabel{found:defn:representable-isofibrations}{{3.4.8}{35}{representably defined isofibrations in 2-categories}{thm.3.4.8}{}}
\newlabel{found:lem:representable-isofibration}{{3.4.9}{36}{}{thm.3.4.9}{}}
\newlabel{found:lem:proj-is-1-conservative}{{3.4.10}{36}{}{thm.3.4.10}{}}
\newlabel{found:eq:equiv-to-lift}{{3.4.11}{36}{}{equation.3.4.11}{}}
\newlabel{found:cor:recog-fibred-equivs}{{3.4.12}{36}{}{thm.3.4.12}{}}
\newlabel{found:defn:fibred-equivalence}{{3.4.13}{37}{fibred equivalence}{thm.3.4.13}{}}
\newlabel{found:obs:1cell-ind-uniqueness-reloaded}{{3.5.1}{37}{uniqueness of 1-cell induction revisited}{thm.3.5.1}{}}
\newlabel{found:eq:induced-1-cell-comparison}{{3.5.2}{37}{uniqueness of 1-cell induction revisited}{equation.3.5.2}{}}
\newlabel{found:obs:squares-set}{{3.5.3}{38}{}{thm.3.5.3}{}}
\newlabel{found:obs:groupoid-components}{{3.5.4}{38}{}{thm.3.5.4}{}}
\newlabel{found:lem:sq-as-a-functor}{{3.5.5}{38}{}{thm.3.5.5}{}}
\newlabel{found:eq:the-real-sq-functor}{{3.5.6}{38}{}{equation.3.5.6}{}}
\newlabel{found:lem:cpts-and-comma-2-cells}{{3.5.7}{39}{}{thm.3.5.7}{}}
\newlabel{found:eq:first-comma-cone}{{3.5.8}{39}{}{equation.3.5.8}{}}
\newlabel{found:eq:other-comma-cone}{{3.5.9}{39}{}{equation.3.5.9}{}}
\newlabel{found:sec:qcatadj}{{4}{40}{Adjunctions of quasi-categories}{section.4}{}}
\newlabel{found:defn:adjunction}{{4.0.1}{40}{adjunction}{thm.4.0.1}{}}
\newlabel{found:ex:simp.quillen.adj}{{4.0.4}{41}{}{thm.4.0.4}{}}
\newlabel{found:subsec:RARI}{{4.1}{42}{Right adjoint right inverse adjunctions}{subsection.4.1}{}}
\newlabel{found:defn:RARI}{{4.1.1}{42}{}{thm.4.1.1}{}}
\newlabel{found:lem:adjunction.from.isos}{{4.1.2}{42}{}{thm.4.1.2}{}}
\newlabel{found:eq:triangle-calculation-1}{{4.1.3}{42}{Right adjoint right inverse adjunctions}{equation.4.1.3}{}}
\newlabel{found:eq:triangle-calculation-2}{{4.1.4}{43}{Right adjoint right inverse adjunctions}{equation.4.1.4}{}}
\newlabel{found:rmk:idempotent-isomorphisms}{{4.1.5}{43}{idempotent isomorphisms}{thm.4.1.5}{}}
\newlabel{found:lem:technicalsliceadjunction}{{4.1.6}{43}{}{thm.4.1.6}{}}
\newlabel{found:eq:technicalsliceadjunction}{{4.1.7}{43}{}{equation.4.1.7}{}}
\newlabel{found:lem:isofibration-RARI}{{4.1.8}{44}{}{thm.4.1.8}{}}
\newlabel{found:subsec:terminal}{{4.2}{44}{Terminal objects as adjoint functors}{subsection.4.2}{}}
\newlabel{found:defn:terminal}{{4.2.1}{44}{terminal objects}{thm.4.2.1}{}}
\newlabel{found:ex:slice-terminal}{{4.2.2}{45}{slices have terminal objects}{thm.4.2.2}{}}
\newlabel{found:lem:min-term-pres}{{4.2.3}{45}{minimal information required to display a terminal object}{thm.4.2.3}{}}
\newlabel{found:lem:adj-ext-univ}{{4.2.4}{45}{}{thm.4.2.4}{}}
\newlabel{found:prop:terminal-ext-univ}{{4.2.5}{46}{}{thm.4.2.5}{}}
\newlabel{found:ex:terminaldefn}{{4.2.6}{46}{}{thm.4.2.6}{}}
\newlabel{found:prop:adj-comp}{{4.3.1}{47}{}{thm.4.3.1}{}}
\newlabel{found:prop:equivtoadjoint}{{4.3.2}{47}{}{thm.4.3.2}{}}
\newlabel{found:prop:expadj}{{4.3.3}{47}{}{thm.4.3.3}{}}
\newlabel{found:prop:terminaldefn}{{4.3.4}{48}{}{thm.4.3.4}{}}
\newlabel{found:eq:commaobjdefn}{{4.4.1}{48}{The universal property of adjunctions}{equation.4.4.1}{}}
\newlabel{found:prop:adjointequiv}{{4.4.2}{48}{}{thm.4.4.2}{}}
\newlabel{found:prop:adjointequivconverse}{{4.4.3}{49}{}{thm.4.4.3}{}}
\newlabel{found:obs:pointwise-adjoint-correspondence}{{4.4.4}{50}{the hom-spaces of a quasi-category}{thm.4.4.4}{}}
\newlabel{found:rmk:vs-lurie-adjunction}{{4.4.5}{50}{}{thm.4.4.5}{}}
\newlabel{found:lem:slice-equiv-comma}{{4.4.6}{50}{}{thm.4.4.6}{}}
\newlabel{found:prop:terminalconverse}{{4.4.7}{51}{}{thm.4.4.7}{}}
\newlabel{found:prop:pointwise-univ-adj}{{4.4.8}{51}{the pointwise universal property of an adjunction}{thm.4.4.8}{}}
\newlabel{found:obs:universal-property-of-epsilon}{{4.4.10}{52}{unpacking this pointwise universal property of an adjunction}{thm.4.4.10}{}}
\newlabel{found:lem:RARI-lifting}{{4.4.12}{53}{right adjoint right inverse as a lifting property}{thm.4.4.12}{}}
\newlabel{found:eq:RARI-lifting}{{4.4.13}{53}{right adjoint right inverse as a lifting property}{equation.4.4.13}{}}
\newlabel{found:subsec:fibred.adjunction}{{4.5}{54}{Fibred adjunctions}{subsection.4.5}{}}
\newlabel{found:defn:fibred.adj}{{4.5.1}{54}{fibred adjunctions}{thm.4.5.1}{}}
\newlabel{found:lem:missed-lemma}{{4.5.2}{55}{}{thm.4.5.2}{}}
\newlabel{found:cor:missed-lemma}{{4.5.3}{55}{}{thm.4.5.3}{}}
\newlabel{found:ex:fibred-technical-slice-adjunction}{{4.5.4}{55}{}{thm.4.5.4}{}}
\newlabel{found:ex:isofib-section.fibred.adjunction}{{4.5.5}{55}{fibred isofibration RARIs}{thm.4.5.5}{}}
\newlabel{found:eq:fibred.terminal}{{4.5.6}{55}{fibred isofibration RARIs}{equation.4.5.6}{}}
\newlabel{found:obs:isofib-section.fibred.adjunction}{{4.5.7}{56}{}{thm.4.5.7}{}}
\newlabel{found:ex:comp.ident.adj}{{4.5.8}{56}{}{thm.4.5.8}{}}
\newlabel{found:eq:comp.ident.adj}{{4.5.9}{57}{}{equation.4.5.9}{}}
\newlabel{found:sec:limits}{{5}{57}{Limits and colimits}{section.5}{}}
\newlabel{found:defn:abs-right-lift}{{5.0.1}{57}{}{thm.5.0.1}{}}
\newlabel{found:eq:absRlifting}{{5.0.2}{57}{}{equation.5.0.2}{}}
\newlabel{found:eq:abs-lifting-property}{{5.0.3}{57}{}{equation.5.0.3}{}}
\newlabel{found:ex:adjasabslifting}{{5.0.4}{58}{}{thm.5.0.4}{}}
\newlabel{found:eq:adjasabslifting}{{5.0.5}{58}{}{equation.5.0.5}{}}
\newlabel{found:eq:unitdefn}{{5.0.6}{58}{Limits and colimits}{equation.5.0.6}{}}
\newlabel{found:eq:comma-cones}{{5.1.1}{58}{Absolute liftings and comma objects}{equation.5.1.1}{}}
\newlabel{found:eq:w-def-prop}{{5.1.2}{59}{Absolute liftings and comma objects}{equation.5.1.2}{}}
\newlabel{found:prop:absliftingtranslation}{{5.1.3}{59}{}{thm.5.1.3}{}}
\newlabel{found:eq:induced-u-from-nattrans-k}{{5.1.4}{60}{Absolute liftings and comma objects}{equation.5.1.4}{}}
\newlabel{found:lem:represented-nat-trans}{{5.1.6}{60}{}{thm.5.1.6}{}}
\newlabel{found:eq:rel-mu-lambda}{{5.1.7}{61}{Absolute liftings and comma objects}{equation.5.1.7}{}}
\newlabel{found:prop:absliftingtranslation2}{{5.1.8}{61}{}{thm.5.1.8}{}}
\newlabel{found:eq:absRlifting.2}{{5.1.9}{62}{}{equation.5.1.9}{}}
\newlabel{found:eq:induced-from-lambda}{{5.1.10}{62}{}{equation.5.1.10}{}}
\newlabel{found:eq:t-for-lim-def-prop}{{5.1.11}{63}{Absolute liftings and comma objects}{equation.5.1.11}{}}
\newlabel{found:prop:right.liftings.as.fibred.terminal.objects}{{5.1.12}{63}{}{thm.5.1.12}{}}
\newlabel{found:eq:fibred.terminal.2}{{5.1.13}{63}{}{equation.5.1.13}{}}
\newlabel{found:eq:limit-lifting-prop}{{5.1.14}{63}{Absolute liftings and comma objects}{equation.5.1.14}{}}
\newlabel{found:eq:m-def-prop}{{5.1.15}{64}{Absolute liftings and comma objects}{equation.5.1.15}{}}
\newlabel{found:obs:right.liftings.as.fibred.terminal.objects}{{5.1.16}{64}{}{thm.5.1.16}{}}
\newlabel{found:prop:translated.lifting}{{5.1.17}{64}{}{thm.5.1.17}{}}
\newlabel{found:eq:orig.lifting}{{5.1.18}{64}{}{equation.5.1.18}{}}
\newlabel{found:eq:translated.lifting}{{5.1.19}{64}{}{equation.5.1.19}{}}
\newlabel{found:eq:ran.as.fibred.adj}{{5.1.20}{65}{Absolute liftings and comma objects}{equation.5.1.20}{}}
\newlabel{found:eq:adj.for.trans}{{5.1.21}{65}{Absolute liftings and comma objects}{equation.5.1.21}{}}
\newlabel{found:defn:limit}{{5.2.2}{66}{}{thm.5.2.2}{}}
\newlabel{found:eq:genericlimit}{{5.2.3}{66}{}{equation.5.2.3}{}}
\newlabel{found:eq:genericcolimit}{{5.2.4}{66}{}{equation.5.2.4}{}}
\newlabel{found:prop:limits.as.terminal.objects}{{5.2.6}{67}{}{thm.5.2.6}{}}
\newlabel{found:lem:cone-equiv-fatcone}{{5.2.7}{67}{}{thm.5.2.7}{}}
\newlabel{found:prop:limits.are.limits}{{5.2.8}{68}{}{thm.5.2.8}{}}
\newlabel{found:defn:families.of.diagrams}{{5.2.9}{68}{}{thm.5.2.9}{}}
\newlabel{found:eq:limits.of.a.family}{{5.2.10}{68}{}{equation.5.2.10}{}}
\newlabel{found:prop:families.of.diagrams}{{5.2.11}{68}{}{thm.5.2.11}{}}
\newlabel{found:prop:limitsasadjunctions}{{5.2.12}{69}{}{thm.5.2.12}{}}
\newlabel{found:prop:RAPL}{{5.2.13}{69}{}{thm.5.2.13}{}}
\newlabel{found:cor:equivprescolim}{{5.2.14}{70}{}{thm.5.2.14}{}}
\newlabel{found:obs:transpose-abs-lifting}{{5.2.15}{70}{}{thm.5.2.15}{}}
\newlabel{found:eq:untransposed}{{5.2.16}{70}{}{equation.5.2.16}{}}
\newlabel{found:eq:transposed}{{5.2.17}{70}{}{equation.5.2.17}{}}
\newlabel{found:prop:pointwise-limits-in-functor-quasi-categories}{{5.2.18}{71}{pointwise limits in functor quasi-categories}{thm.5.2.18}{}}
\newlabel{found:prop:ran.adj.limits}{{5.2.19}{71}{}{thm.5.2.19}{}}
\newlabel{found:cor:ran.adj.limits}{{5.2.20}{72}{}{thm.5.2.20}{}}
\newlabel{found:ntn:pb.po.joins}{{5.2.22}{73}{pushout and pullback diagrams}{thm.5.2.22}{}}
\newlabel{found:obs:loops.diag.fam}{{5.2.24}{73}{}{thm.5.2.24}{}}
\newlabel{found:defn:loop.susp}{{5.2.25}{73}{loop spaces and suspensions}{thm.5.2.25}{}}
\newlabel{found:prop:loops-suspension}{{5.2.27}{74}{}{thm.5.2.27}{}}
\newlabel{found:eq:pullback.pushout.adj}{{5.2.28}{74}{Limits and colimits as absolute lifting diagrams}{equation.5.2.28}{}}
\newlabel{found:eq:loop.susp.var}{{5.2.29}{74}{Limits and colimits as absolute lifting diagrams}{equation.5.2.29}{}}
\newlabel{found:thm:splitgeorealizations}{{5.3.1}{76}{}{thm.5.3.1}{}}
\newlabel{found:lem:doms2catlemma}{{5.3.2}{76}{}{thm.5.3.2}{}}
\newlabel{found:sec:pointwise}{{6}{78}{Pointwise universal properties}{section.6}{}}
\newlabel{found:defn:pointwise-abs-lifting}{{6.1.1}{78}{pointwise universal property of absoluting lifting diagrams}{thm.6.1.1}{}}
\newlabel{found:eq:pointwise-absRlifting}{{6.1.2}{78}{pointwise universal property of absoluting lifting diagrams}{equation.6.1.2}{}}
\newlabel{found:lem:pointwise-terminal}{{6.1.3}{78}{}{thm.6.1.3}{}}
\newlabel{found:thm:pointwise}{{6.1.4}{79}{}{thm.6.1.4}{}}
\newlabel{found:eq:term.obj.assumption}{{6.1.5}{79}{Pointwise absolute lifting}{equation.6.1.5}{}}
\newlabel{found:eq:term.obj.extension}{{6.1.6}{79}{Pointwise absolute lifting}{equation.6.1.6}{}}
\newlabel{found:eq:term.obj.lifting}{{6.1.7}{79}{Pointwise absolute lifting}{equation.6.1.7}{}}
\newlabel{found:cor:pointwise}{{6.1.8}{80}{}{thm.6.1.8}{}}
\newlabel{found:thm:simplicial-Quillen-adjunction}{{6.2.1}{81}{}{thm.6.2.1}{}}
\newlabel{found:eq:simp.adj.extension}{{6.2.2}{82}{Simplicial Quillen adjunctions are adjunctions of quasi-categories}{equation.6.2.2}{}}
\newlabel{found:tocindent-1}{0pt}
\newlabel{found:tocindent0}{68.84703pt}
\newlabel{found:tocindent1}{77.65321pt}
\newlabel{found:tocindent2}{37.93518pt}
\newlabel{found:tocindent3}{0pt}

%% file: cohadj_lab.tex
\newlabel{cohadj:sec:intro-adj-data}{{1.1}{2}{Adjunction data}{subsection.1.1}{}}
\newlabel{cohadj:eq:sampletriangleidentities}{{1.1.1}{3}{Adjunction data}{equation.1.1.1}{}}
\newlabel{cohadj:eq:horn-1}{{1.1.2}{3}{Adjunction data}{equation.1.1.2}{}}
\newlabel{cohadj:eq:middle.four}{{1.1.3}{4}{Adjunction data}{equation.1.1.3}{}}
\newlabel{cohadj:eq:horn-2}{{1.1.4}{4}{Adjunction data}{equation.1.1.4}{}}
\newlabel{cohadj:sec:computads}{{2}{7}{Simplicial computads}{section.2}{}}
\newlabel{cohadj:subsec:computads}{{2.1}{7}{Simplicial categories and simplicial computads}{subsection.2.1}{}}
\newlabel{cohadj:ntn:whiskering}{{2.1.2}{8}{whiskering in a simplicial category}{thm.2.1.2}{}}
\newlabel{cohadj:ntn:generic-cattwo}{{2.1.3}{8}{the generic $n$-arrow}{thm.2.1.3}{}}
\newlabel{cohadj:defn:computads}{{2.1.4}{8}{(relative) simplicial computads}{thm.2.1.4}{}}
\newlabel{cohadj:obs:simp-computad-char}{{2.1.5}{8}{an explicit characterisation of simplicial computads}{thm.2.1.5}{}}
\newlabel{cohadj:eq:computad-arrrow-decomp}{{2.1.6}{9}{an explicit characterisation of simplicial computads}{equation.2.1.6}{}}
\newlabel{cohadj:ex:gothic-C}{{2.1.10}{9}{}{thm.2.1.10}{}}
\newlabel{cohadj:ex:homotopy-coherent-simplex}{{2.1.11}{10}{}{thm.2.1.11}{}}
\newlabel{cohadj:subsec:subcomputads}{{2.2}{10}{Simplicial subcomputads}{subsection.2.2}{}}
\newlabel{cohadj:ex:mono-subcomputad}{{2.2.3}{10}{}{thm.2.2.3}{}}
\newlabel{cohadj:defn:generated-subcomputad}{{2.2.4}{11}{}{thm.2.2.4}{}}
\newlabel{cohadj:ex:simp-cat.(co)skeleta}{{2.2.5}{11}{(co)skeleta of simplicial categories}{thm.2.2.5}{}}
\newlabel{cohadj:prop:simp-computad-maps}{{2.2.6}{11}{}{thm.2.2.6}{}}
\newlabel{cohadj:sec:generic-adj}{{3}{12}{The generic adjunction}{section.3}{}}
\newlabel{cohadj:ssec:graphical}{{3.1}{12}{A graphical calculus for the simplicial category \protect \texorpdfstring {$\protect \Adj $}{Adj}}{subsection.3.1}{}}
\newlabel{cohadj:eq:samplesquiggle}{{3.1.1}{13}{A graphical calculus for the simplicial category \protect \texorpdfstring {$\protect \Adj $}{Adj}}{equation.3.1.1}{}}
\newlabel{cohadj:def:strict-und-squiggles}{{3.1.2}{14}{strictly undulating squiggles}{thm.3.1.2}{}}
\newlabel{cohadj:item:strict-und-squiggle-1}{{{{(i)}}}{14}{strictly undulating squiggles}{Item.1}{}}
\newlabel{cohadj:item:strict-und-squiggle-2}{{{{(ii)}}}{14}{strictly undulating squiggles}{Item.2}{}}
\newlabel{cohadj:item:und-squiggle-2}{{{{(ii)$'$}}}{14}{strictly undulating squiggles}{Item.3}{}}
\newlabel{cohadj:def:comp-squiggles}{{3.1.3}{14}{composing squiggles}{thm.3.1.3}{}}
\newlabel{cohadj:eq:samplesquiggle2}{{3.1.6}{15}{simplicial action on strictly undulating squiggles}{equation.3.1.6}{}}
\newlabel{cohadj:obs:interval.rep}{{3.1.7}{17}{interval representation}{thm.3.1.7}{}}
\newlabel{cohadj:obs:squiggle-vertices}{{3.1.9}{18}{the vertices of an arrow in $\Adj $}{thm.3.1.9}{}}
\newlabel{cohadj:prop:adjcomputad}{{3.1.10}{19}{}{thm.3.1.10}{}}
\newlabel{cohadj:ex:sample-adjunction-data}{{3.1.11}{19}{adjunction data in $\Adj $}{thm.3.1.11}{}}
\newlabel{cohadj:eq:ualpha-faces}{{3.1.12}{20}{adjunction data in $\Adj $}{equation.3.1.12}{}}
\newlabel{cohadj:ssec:2-cat-adj}{{3.2}{21}{The simplicial category \texorpdfstring {$\protect \Adj $}{Adj} as a 2-category}{subsection.3.2}{}}
\newlabel{cohadj:prop:adj.2-cat}{{3.2.2}{21}{}{thm.3.2.2}{}}
\newlabel{cohadj:eq:samplesquiggles4}{{3.2.3}{21}{The simplicial category \texorpdfstring {$\protect \Adj $}{Adj} as a 2-category}{equation.3.2.3}{}}
\newlabel{cohadj:eq:samplesquiggle3}{{3.2.4}{22}{The simplicial category \texorpdfstring {$\protect \Adj $}{Adj} as a 2-category}{equation.3.2.4}{}}
\newlabel{cohadj:rmk:hammock}{{3.2.5}{22}{}{thm.3.2.5}{}}
\newlabel{cohadj:eq:karol}{{3.2.6}{23}{}{equation.3.2.6}{}}
\newlabel{cohadj:ssec:adj-is-adj}{{3.3}{23}{The 2-categorical universal property of \protect \texorpdfstring {$\protect \Adj $}{Adj}}{subsection.3.3}{}}
\newlabel{cohadj:obs:2-cat.as.simp-cat}{{3.3.1}{24}{2-categories as simplicial categories}{thm.3.3.1}{}}
\newlabel{cohadj:eq:genadjdata}{{3.3.3}{24}{the adjunction in \protect \texorpdfstring {$\protect \Adj $}{Adj}}{equation.3.3.3}{}}
\newlabel{cohadj:prop:2-cat-univ-Adj}{{3.3.4}{24}{a 2-categorical universal property of \protect \texorpdfstring {$\protect \Adj $}{Adj}}{thm.3.3.4}{}}
\newlabel{cohadj:cor:schanuel-street-iso}{{3.3.5}{25}{}{thm.3.3.5}{}}
\newlabel{cohadj:obs:delta.adj.duality}{{3.3.6}{26}{adjunctions in \protect \texorpdfstring {$\protect \Del +$}{Delta+}}{thm.3.3.6}{}}
\newlabel{cohadj:eq:elem.op.adj}{{3.3.7}{26}{adjunctions in \protect \texorpdfstring {$\protect \Del +$}{Delta+}}{equation.3.3.7}{}}
\newlabel{cohadj:rmk:schanuel-street}{{3.3.8}{26}{the Schanuel and Street 2-category $\Adj $}{thm.3.3.8}{}}
\newlabel{cohadj:sec:adjunction-data}{{4}{28}{Adjunction data}{section.4}{}}
\newlabel{cohadj:ssec:fillable}{{4.1}{28}{Fillable arrows}{subsection.4.1}{}}
\newlabel{cohadj:obs:fillable-dist-face}{{4.1.2}{29}{fillable arrows and distinguished faces}{thm.4.1.2}{}}
\newlabel{cohadj:lem:fillablecodim1}{{4.1.3}{29}{}{thm.4.1.3}{}}
\newlabel{cohadj:lem:fillablecodim1'}{{4.1.4}{30}{}{thm.4.1.4}{}}
\newlabel{cohadj:ssec:parental}{{4.2}{30}{Parental subcomputads}{subsection.4.2}{}}
\newlabel{cohadj:defn:parental-subcomputad}{{4.2.2}{30}{parental subcomputads of $\Adj $}{thm.4.2.2}{}}
\newlabel{cohadj:ex:parental-subcomputad}{{4.2.3}{31}{}{thm.4.2.3}{}}
\newlabel{cohadj:ex:non-parental-subcomputad}{{4.2.4}{31}{a non-example}{thm.4.2.4}{}}
\newlabel{cohadj:ex:another-parental-subcomputad}{{4.2.5}{31}{}{thm.4.2.5}{}}
\newlabel{cohadj:ntn:generic-catthree}{{4.2.6}{31}{}{thm.4.2.6}{}}
\newlabel{cohadj:eq:3[i]-square}{{4.2.7}{32}{}{equation.4.2.7}{}}
\newlabel{cohadj:def:funct-for-fillable}{{4.2.9}{32}{}{thm.4.2.9}{}}
\newlabel{cohadj:lem:ext-par-subcomp}{{4.2.10}{32}{extending parental subcomputads}{thm.4.2.10}{}}
\newlabel{cohadj:eq:parental-adj-pushout-1}{{4.2.11}{33}{extending parental subcomputads}{equation.4.2.11}{}}
\newlabel{cohadj:eq:parental-adj-pushout-2}{{4.2.12}{33}{extending parental subcomputads}{equation.4.2.12}{}}
\newlabel{cohadj:cor:ext-par-subcomp}{{4.2.13}{34}{}{thm.4.2.13}{}}
\newlabel{cohadj:eq:big-pushout}{{4.2.14}{34}{}{equation.4.2.14}{}}
\newlabel{cohadj:prop:ext-par-subcomp}{{4.2.15}{34}{}{thm.4.2.15}{}}
\newlabel{cohadj:itm:one}{{{{(i)}}}{34}{}{Item.4}{}}
\newlabel{cohadj:itm:two}{{{{(ii)}}}{34}{}{Item.5}{}}
\newlabel{cohadj:ssec:hocoh}{{4.3}{36}{Homotopy Coherent Adjunctions}{subsection.4.3}{}}
\newlabel{cohadj:obs:adjunctions-in-qcat-cats}{{4.3.2}{36}{adjunctions in a quasi-categorically enriched category}{thm.4.3.2}{}}
\newlabel{cohadj:obs:int-univ}{{4.3.4}{37}{the internal universal property of the counit}{thm.4.3.4}{}}
\newlabel{cohadj:lem:rel-lift-terminal}{{4.3.5}{38}{a relative universal property of terminal objects}{thm.4.3.5}{}}
\newlabel{cohadj:prop:rel-int-univ}{{4.3.6}{39}{the relative internal universal property of the counit}{thm.4.3.6}{}}
\newlabel{cohadj:obs:ext-quasicat-homs}{{4.3.7}{39}{}{thm.4.3.7}{}}
\newlabel{cohadj:thm:hty-coherence-exist}{{4.3.8}{40}{}{thm.4.3.8}{}}
\newlabel{cohadj:thm:hty-coherence-exist-I}{{4.3.9}{41}{homotopy coherence of adjunctions I}{thm.4.3.9}{}}
\newlabel{cohadj:thm:hty-coherence-exist-II}{{4.3.11}{42}{homotopy coherence of adjunctions II}{thm.4.3.11}{}}
\newlabel{cohadj:def:adjunction-data}{{4.3.13}{42}{}{thm.4.3.13}{}}
\newlabel{cohadj:ssec:uniqueness}{{4.4}{42}{Homotopical uniqueness of homotopy coherent adjunctions}{subsection.4.4}{}}
\newlabel{cohadj:obs:icon-defn}{{4.4.1}{43}{simplicial enrichment of simplicial categories}{thm.4.4.1}{}}
\newlabel{cohadj:lem:isofib-icon}{{4.4.2}{43}{}{thm.4.4.2}{}}
\newlabel{cohadj:lem:conservative-icon}{{4.4.4}{45}{}{thm.4.4.4}{}}
\newlabel{cohadj:lem:cohadj.space.Kan}{{4.4.6}{46}{}{thm.4.4.6}{}}
\newlabel{cohadj:prop:hty-uniqueness}{{4.4.7}{46}{}{thm.4.4.7}{}}
\newlabel{cohadj:eq:trans-lift-prob-1}{{4.4.8}{47}{Homotopical uniqueness of homotopy coherent adjunctions}{equation.4.4.8}{}}
\newlabel{cohadj:obs:epsilon-icon-expl}{{4.4.9}{47}{}{thm.4.4.9}{}}
\newlabel{cohadj:thm:hty-uniqueness-I}{{4.4.11}{48}{}{thm.4.4.11}{}}
\newlabel{cohadj:prop:counits-proj-trivial}{{4.4.12}{48}{}{thm.4.4.12}{}}
\newlabel{cohadj:lem:fib-contract-fibres}{{4.4.13}{48}{}{thm.4.4.13}{}}
\newlabel{cohadj:prop:leftadjs-proj-trivial}{{4.4.17}{49}{}{thm.4.4.17}{}}
\newlabel{cohadj:thm:hty-uniqueness-II}{{4.4.18}{50}{}{thm.4.4.18}{}}
\newlabel{cohadj:sec:weighted}{{5}{50}{Weighted limits in \protect \texorpdfstring {$\protect \qCat _\infty $}{qCat}}{section.5}{}}
\newlabel{cohadj:subsec:weighted}{{5.1}{50}{Weighted limits and colimits}{subsection.5.1}{}}
\newlabel{cohadj:defn:cotensors}{{5.1.1}{50}{cotensors}{thm.5.1.1}{}}
\newlabel{cohadj:eq:wlimformula}{{5.1.3}{51}{weighted limits}{equation.5.1.3}{}}
\newlabel{cohadj:eq:wlim-def-prop}{{5.1.4}{51}{weighted limits}{equation.5.1.4}{}}
\newlabel{cohadj:ex:rep-weights}{{5.1.5}{51}{representable weights}{thm.5.1.5}{}}
\newlabel{cohadj:eq:yoneda}{{5.1.6}{51}{representable weights}{equation.5.1.6}{}}
\newlabel{cohadj:obs:weighted-cocontinuity}{{5.1.7}{51}{}{thm.5.1.7}{}}
\newlabel{cohadj:ex:diagrams-weighted-limit}{{5.1.8}{52}{diagrams}{thm.5.1.8}{}}
\newlabel{cohadj:ex:commaweightedlimit}{{5.1.9}{52}{comma quasi-categories}{thm.5.1.9}{}}
\newlabel{cohadj:ex:homotopyweighted}{{5.1.10}{52}{homotopy limits as weighted limits}{thm.5.1.10}{}}
\newlabel{cohadj:lem:lanweights}{{5.1.11}{52}{weighted limits and Kan extensions}{thm.5.1.11}{}}
\newlabel{cohadj:subsec:weighted-qcat}{{5.2}{53}{Weighted limits in the quasi-categorical context}{subsection.5.2}{}}
\newlabel{cohadj:defn:proj-cof}{{5.2.1}{53}{projective cofibrations}{thm.5.2.1}{}}
\newlabel{cohadj:prop:projwlims2}{{5.2.2}{53}{}{thm.5.2.2}{}}
\newlabel{cohadj:eq:layer-in-tower}{{5.2.3}{53}{Weighted limits in the quasi-categorical context}{equation.5.2.3}{}}
\newlabel{cohadj:prop:projwlims}{{5.2.4}{54}{}{thm.5.2.4}{}}
\newlabel{cohadj:prop:proj-wlim-homotopical}{{5.2.6}{54}{}{thm.5.2.6}{}}
\newlabel{cohadj:rmk:projwlims}{{5.2.7}{54}{}{thm.5.2.7}{}}
\newlabel{cohadj:rmk:categories-special-case}{{5.2.9}{55}{2-categorical weighted limits and quasi-categorical weighted limits}{thm.5.2.9}{}}
\newlabel{cohadj:subsec:collage}{{5.3}{55}{The collage construction}{subsection.5.3}{}}
\newlabel{cohadj:obs:coll-right-adj}{{5.3.2}{55}{a right adjoint to the collage construction}{thm.5.3.2}{}}
\newlabel{cohadj:prop:projcofchar2}{{5.3.3}{56}{}{thm.5.3.3}{}}
\newlabel{cohadj:eq:pushout-transform}{{5.3.4}{56}{The collage construction}{equation.5.3.4}{}}
\newlabel{cohadj:prop:projcofchar}{{5.3.5}{58}{}{thm.5.3.5}{}}
\newlabel{cohadj:sec:formal}{{6}{58}{The formal theory of homotopy coherent monads}{section.6}{}}
\newlabel{cohadj:ssec:formalmonads}{{6.1}{59}{Weighted limits for the formal theory of monads}{subsection.6.1}{}}
\newlabel{cohadj:obs:mnd-weights}{{6.1.3}{59}{weights on $\Mnd $}{thm.6.1.3}{}}
\newlabel{cohadj:defn:W+}{{6.1.4}{60}{monad resolutions}{thm.6.1.4}{}}
\newlabel{cohadj:eq:resolution2}{{6.1.5}{60}{monad resolutions}{equation.6.1.5}{}}
\newlabel{cohadj:defn:W-}{{6.1.6}{61}{}{thm.6.1.6}{}}
\newlabel{cohadj:defn:EMobject}{{6.1.7}{61}{quasi-category of algebras}{thm.6.1.7}{}}
\newlabel{cohadj:lem:projcof1}{{6.1.8}{61}{}{thm.6.1.8}{}}
\newlabel{cohadj:rmk:expl-htycoh-alg}{{6.1.10}{61}{}{thm.6.1.10}{}}
\newlabel{cohadj:eq:expl-htycoh-alg-cond}{{6.1.11}{61}{}{equation.6.1.11}{}}
\newlabel{cohadj:eq:algebraresolution}{{6.1.12}{62}{}{equation.6.1.12}{}}
\newlabel{cohadj:ex:monadicadj}{{6.1.14}{63}{monadic adjunction}{thm.6.1.14}{}}
\newlabel{cohadj:eq:monadicadjunctionweights}{{6.1.15}{63}{monadic adjunction}{equation.6.1.15}{}}
\newlabel{cohadj:ssec:monadicequivs}{{6.2}{63}{Conservativity of the monadic forgetful functor}{subsection.6.2}{}}
\newlabel{cohadj:prop:conservative-char}{{6.2.2}{63}{}{thm.6.2.2}{}}
\newlabel{cohadj:cor:uTconservative}{{6.2.3}{64}{}{thm.6.2.3}{}}
\newlabel{cohadj:obs:understanding-monadic-forgetful-functor}{{6.2.4}{64}{}{thm.6.2.4}{}}
\newlabel{cohadj:rmk:u-map}{{6.2.5}{65}{}{thm.6.2.5}{}}
\newlabel{cohadj:ssec:algcolims}{{6.3}{65}{Colimit representation of algebras}{subsection.6.3}{}}
\newlabel{cohadj:eq:splitcoeq}{{6.3.1}{65}{Colimit representation of algebras}{equation.6.3.1}{}}
\newlabel{cohadj:eq:splitgeorealizations}{{6.3.2}{66}{}{equation.6.3.2}{}}
\newlabel{cohadj:rec:splitgeorealizations}{{6.3.3}{66}{constructing the triangle in theorem \refI {thm:splitgeorealizations}}{thm.6.3.3}{}}
\newlabel{cohadj:eq:canonical-2-cell}{{6.3.4}{66}{constructing the triangle in theorem \refI {thm:splitgeorealizations}}{equation.6.3.4}{}}
\newlabel{cohadj:defn:u-split-aug}{{6.3.5}{66}{$u$-split augmented simplicial objects}{thm.6.3.5}{}}
\newlabel{cohadj:eq:u-split-pullback}{{6.3.6}{67}{$u$-split augmented simplicial objects}{equation.6.3.6}{}}
\newlabel{cohadj:prop:uTcreates}{{6.3.7}{67}{}{thm.6.3.7}{}}
\newlabel{cohadj:thm:uTcreates}{{6.3.8}{67}{}{thm.6.3.8}{}}
\newlabel{cohadj:obs:uTcreates}{{6.3.9}{67}{}{thm.6.3.9}{}}
\newlabel{cohadj:eq:colim.diag.1}{{6.3.10}{67}{}{equation.6.3.10}{}}
\newlabel{cohadj:eq:colimit-diag-ut-split}{{6.3.11}{69}{Colimit representation of algebras}{equation.6.3.11}{}}
\newlabel{cohadj:eq:colimit-diag-ut-split-under}{{6.3.12}{69}{Colimit representation of algebras}{equation.6.3.12}{}}
\newlabel{cohadj:eq:Wscdefn}{{6.3.14}{70}{a direct description of $S(u^t)$}{equation.6.3.14}{}}
\newlabel{cohadj:obs:canonical.alg.pres}{{6.3.15}{70}{}{thm.6.3.15}{}}
\newlabel{cohadj:eq:usplitpushout}{{6.3.16}{70}{}{equation.6.3.16}{}}
\newlabel{cohadj:thm:colim-rep-algebras}{{6.3.17}{71}{canonical colimit representation of algebras}{thm.6.3.17}{}}
\newlabel{cohadj:eq:canonical-colimits}{{6.3.18}{71}{canonical colimit representation of algebras}{equation.6.3.18}{}}
\newlabel{cohadj:eq:downstairs-colimit}{{6.3.19}{71}{canonical colimit representation of algebras}{equation.6.3.19}{}}
\newlabel{cohadj:sec:monadicity}{{7}{72}{Monadicity}{section.7}{}}
\newlabel{cohadj:ssec:comparison}{{7.1}{72}{Comparison with the monadic adjunction}{subsection.7.1}{}}
\newlabel{cohadj:lem:actionquotient}{{7.1.3}{73}{}{thm.7.1.3}{}}
\newlabel{cohadj:eq:actionquotient}{{7.1.4}{73}{}{equation.7.1.4}{}}
\newlabel{cohadj:obs:lanW-}{{7.1.5}{73}{}{thm.7.1.5}{}}
\newlabel{cohadj:eq:monadicadjunctionweight}{{7.1.6}{74}{}{equation.7.1.6}{}}
\newlabel{cohadj:eq:monadiccomparison}{{7.1.10}{74}{comparison with the monadic adjunction}{equation.7.1.10}{}}
\newlabel{cohadj:ssec:monadicity}{{7.2}{74}{The monadicity theorem}{subsection.7.2}{}}
\newlabel{cohadj:eq:usplit}{{7.2.2}{75}{$u$-split simplicial objects}{equation.7.2.2}{}}
\newlabel{cohadj:eq:usplithyp}{{7.2.3}{75}{$u$-split simplicial objects}{equation.7.2.3}{}}
\newlabel{cohadj:thm:monadiccomparisonadj}{{7.2.4}{75}{monadicity I}{thm.7.2.4}{}}
\newlabel{cohadj:eq:Ldefn}{{7.2.5}{75}{The monadicity theorem}{equation.7.2.5}{}}
\newlabel{cohadj:eq:Labs-lifting}{{7.2.6}{76}{The monadicity theorem}{equation.7.2.6}{}}
\newlabel{cohadj:thm:monadicity}{{7.2.7}{77}{monadicity II}{thm.7.2.7}{}}
\newlabel{cohadj:tocindent-1}{0pt}
\newlabel{cohadj:tocindent0}{15.01021pt}
\newlabel{cohadj:tocindent1}{20.88377pt}
\newlabel{cohadj:tocindent2}{34.49158pt}
\newlabel{cohadj:tocindent3}{0pt}

%% file: reedy_lab.tex
\newlabel{reedy:eq:reedy's.po}{{1.1}{2}{Introduction}{equation.1.1}{}}
\newlabel{reedy:eq:reedy's.delta}{{1.2}{3}{Introduction}{equation.1.2}{}}
\newlabel{reedy:ntn:index-stuff}{{1.4}{4}{general index notation}{thm.1.4}{}}
\newlabel{reedy:ex:tensor-cotensor}{{1.6}{5}{tensors and cotensors}{thm.1.6}{}}
\newlabel{reedy:defn:ends-and-coends}{{1.7}{6}{ends and coends}{thm.1.7}{}}
\newlabel{reedy:defn:cat-of-elts}{{1.8}{6}{categories of elements}{thm.1.8}{}}
\newlabel{reedy:obs:coend-in-sets}{{1.9}{7}{coends in sets}{thm.1.9}{}}
\newlabel{reedy:eq:pb-coends-in-sets}{{1.10}{7}{coends in sets}{equation.1.10}{}}
\newlabel{reedy:eq:wlimformula}{{1.12}{7}{weighted limits and colimits}{equation.1.12}{}}
\newlabel{reedy:ex:weighted-conical}{{1.13}{7}{terminal weights}{thm.1.13}{}}
\newlabel{reedy:ex:weighted-yoneda}{{1.14}{8}{representable weights}{thm.1.14}{}}
\newlabel{reedy:eq:yoneda}{{1.15}{8}{representable weights}{equation.1.15}{}}
\newlabel{reedy:ex:matching.preview}{{1.17}{8}{}{thm.1.17}{}}
\newlabel{reedy:obs:weighted.as.ordinary}{{1.19}{8}{weighted limits as ordinary limits}{thm.1.19}{}}
\newlabel{reedy:obs:wcolim-in-sets}{{1.21}{8}{weighted colimits in sets}{thm.1.21}{}}
\newlabel{reedy:sec:reedy}{{2}{9}{Reedy categories}{section.2}{}}
\newlabel{reedy:defn:reedy}{{2.1}{9}{Reedy categories}{thm.2.1}{}}
\newlabel{reedy:ex:poset.reedy}{{2.3}{9}{finite posets}{thm.2.3}{}}
\newlabel{reedy:ex:pushout.reedy.alt}{{2.4}{9}{the pushout diagram}{thm.2.4}{}}
\newlabel{reedy:ex:parallel.pair.reedy}{{2.5}{10}{the parallel pair}{thm.2.5}{}}
\newlabel{reedy:defn:factorisations}{{2.8}{10}{categories of factorisations}{thm.2.8}{}}
\newlabel{reedy:lem:reedy-fact-connect}{{2.9}{10}{}{thm.2.9}{}}
\newlabel{reedy:sec:latching}{{3}{11}{Latching and matching objects}{section.3}{}}
\newlabel{reedy:rec:skel-coskel}{{3.1}{11}{skeleta and coskeleta}{thm.3.1}{}}
\newlabel{reedy:eq:4}{{3.2}{12}{skeleta and coskeleta}{equation.3.2}{}}
\newlabel{reedy:obs:skel-colim}{{3.3}{12}{}{thm.3.3}{}}
\newlabel{reedy:eq:5}{{3.4}{12}{}{equation.3.4}{}}
\newlabel{reedy:lem:skeleta-colimit}{{3.5}{13}{}{thm.3.5}{}}
\newlabel{reedy:eq:6}{{3.6}{13}{}{equation.3.6}{}}
\newlabel{reedy:obs:skel-reps}{{3.7}{13}{skeleta of representables and factorisations}{thm.3.7}{}}
\newlabel{reedy:eq:sk-rep-coend}{{3.8}{13}{skeleta of representables and factorisations}{equation.3.8}{}}
\newlabel{reedy:lem:inductive.defn}{{3.10}{15}{inductive definition of diagrams}{thm.3.10}{}}
\newlabel{reedy:obs:inductive.defn}{{3.11}{15}{inductive definition of natural transformations}{thm.3.11}{}}
\newlabel{reedy:ex:sequence-induction}{{3.12}{15}{}{thm.3.12}{}}
\newlabel{reedy:defn:latching}{{3.13}{16}{latching and matching objects}{thm.3.13}{}}
\newlabel{reedy:ex:simp-latching-1}{{3.14}{16}{}{thm.3.14}{}}
\newlabel{reedy:obs:weights-latching}{{3.15}{16}{weights for latching and matching objects}{thm.3.15}{}}
\newlabel{reedy:eq:7}{{3.16}{16}{weights for latching and matching objects}{equation.3.16}{}}
\newlabel{reedy:lem:sk-rep-reedy}{{3.17}{17}{skeleta of the representables of a Reedy category}{thm.3.17}{}}
\newlabel{reedy:obs:cons-skeleta}{{3.18}{17}{characterising the boundary of a representable}{thm.3.18}{}}
\newlabel{reedy:ex:sequence.latching}{{3.19}{17}{}{thm.3.19}{}}
\newlabel{reedy:ex:pushout.latching}{{3.20}{17}{}{thm.3.20}{}}
\newlabel{reedy:ex:parallel.pair.latching}{{3.21}{18}{}{thm.3.21}{}}
\newlabel{reedy:ex:simp-latching-2}{{3.22}{18}{}{thm.3.22}{}}
\newlabel{reedy:obs:latching.ordinary.colimit}{{3.23}{18}{latching objects as ordinary colimits}{thm.3.23}{}}
\newlabel{reedy:sec:Leibniz-Reedy}{{4}{19}{Leibniz constructions and the Reedy model structure}{section.4}{}}
\newlabel{reedy:obs:box-product}{{4.2}{19}{Leibniz's formula}{thm.4.2}{}}
\newlabel{reedy:defn:leibniz}{{4.4}{20}{the Leibniz construction}{thm.4.4}{}}
\newlabel{reedy:eq:13}{{4.5}{20}{the Leibniz construction}{equation.4.5}{}}
\newlabel{reedy:ex:boundary-prod}{{4.6}{20}{}{thm.4.6}{}}
\newlabel{reedy:obs:leibniz-isos}{{4.7}{20}{Leibniz and structural isomorphisms}{thm.4.7}{}}
\newlabel{reedy:lem:leibniz-cocts}{{4.8}{21}{Leibniz and colimit preservation}{thm.4.8}{}}
\newlabel{reedy:eq:2varadj}{{4.9}{21}{\S \ Leibniz constructions}{equation.4.9}{}}
\newlabel{reedy:lem:leibniz-close}{{4.10}{21}{Leibniz and closures}{thm.4.10}{}}
\newlabel{reedy:obs:leibniz-lifting-properties}{{4.11}{22}{Leibniz and lifting properties}{thm.4.11}{}}
\newlabel{reedy:defn:relative-maps}{{4.14}{22}{relative latching and matching maps}{thm.4.14}{}}
\newlabel{reedy:obs:relative.lifting}{{4.16}{23}{relative latching and matching maps and lifting problems}{thm.4.16}{}}
\newlabel{reedy:eq:i-lift-p}{{4.17}{23}{relative latching and matching maps and lifting problems}{equation.4.17}{}}
\newlabel{reedy:thm:model-structure}{{4.18}{24}{the Reedy model structure}{thm.4.18}{}}
\newlabel{reedy:sec:Leibniz-rel-cell-cx}{{5}{24}{Leibniz constructions and cell complexes}{section.5}{}}
\newlabel{reedy:obs:leibniz-comp}{{5.1}{24}{Leibniz and composites}{thm.5.1}{}}
\newlabel{reedy:defn:transfinite-composites}{{5.3}{25}{transfinite composites}{thm.5.3}{}}
\newlabel{reedy:def:rel-cell}{{5.4}{25}{cell complexes}{thm.5.4}{}}
\newlabel{reedy:eq:10}{{5.5}{26}{cell complexes}{equation.5.5}{}}
\newlabel{reedy:ntn:phi-map}{{5.6}{26}{}{thm.5.6}{}}
\newlabel{reedy:lem:leibniz-tcofp}{{5.7}{26}{Leibniz bifunctors and cell complexes}{thm.5.7}{}}
\newlabel{reedy:eq:11}{{5.8}{27}{Leibniz constructions and cell complexes}{equation.5.8}{}}
\newlabel{reedy:obs:leibniz-tcofp}{{5.9}{27}{}{thm.5.9}{}}
\newlabel{reedy:ex:latch-comp}{{5.10}{28}{relative latching maps of composites II}{thm.5.10}{}}
\newlabel{reedy:cor:leibniz-tcofp}{{5.11}{28}{}{thm.5.11}{}}
\newlabel{reedy:prop:leibniz-tcofp}{{5.12}{28}{}{thm.5.12}{}}
\newlabel{reedy:sec:cellular}{{6}{29}{Cellular presentations and Reedy categories}{section.6}{}}
\newlabel{reedy:obs:two-sided}{{6.1}{29}{skeleta of two-sided representables}{thm.6.1}{}}
\newlabel{reedy:obs:building-skeleta-rep}{{6.2}{29}{building up for skeleta of representables}{thm.6.2}{}}
\newlabel{reedy:prop:building-up}{{6.3}{30}{general building up}{thm.6.3}{}}
\newlabel{reedy:eq:skel-seq}{{6.4}{31}{general building up}{equation.6.4}{}}
\newlabel{reedy:eq:c-cell}{{6.5}{31}{general building up}{equation.6.5}{}}
\newlabel{reedy:eq:skel-seq-cells}{{6.6}{31}{Cellular presentations and Reedy categories}{equation.6.6}{}}
\newlabel{reedy:cor:building-up}{{6.7}{32}{}{thm.6.7}{}}
\newlabel{reedy:cor:B-cell-complex}{{6.8}{32}{}{thm.6.8}{}}
\newlabel{reedy:ex:EZ.simp.set}{{6.9}{33}{}{thm.6.9}{}}
\newlabel{reedy:sec:proof}{{7}{33}{Proof of the Reedy model structure}{section.7}{}}
\newlabel{reedy:lem:triv-cof-char}{{7.1}{33}{}{thm.7.1}{}}
\newlabel{reedy:lem:lifting}{{7.3}{34}{lifting}{thm.7.3}{}}
\newlabel{reedy:lem:factorisation}{{7.4}{34}{factorisation}{thm.7.4}{}}
\newlabel{reedy:eq:factorisation-defn}{{7.5}{35}{Proof of the Reedy model structure}{equation.7.5}{}}
\newlabel{reedy:prop:reedy.cof.gen}{{7.7}{35}{}{thm.7.7}{}}
\newlabel{reedy:sec:hocolim}{{8}{36}{Homotopy limits and colimits}{section.8}{}}
\newlabel{reedy:ex:hocoeq}{{8.2}{37}{homotopy coequalisers}{thm.8.2}{}}
\newlabel{reedy:ex:hoeq}{{8.3}{37}{homotopy equalisers}{thm.8.3}{}}
\newlabel{reedy:prop:reedy.model.dual}{{8.4}{38}{}{thm.8.4}{}}
\newlabel{reedy:ex:mapping}{{8.5}{38}{mapping telescopes}{thm.8.5}{}}
\newlabel{reedy:eq:sequence}{{8.6}{38}{mapping telescopes}{equation.8.6}{}}
\newlabel{reedy:eq:sequence.replacement}{{8.7}{38}{mapping telescopes}{equation.8.7}{}}
\newlabel{reedy:ex:hopushout}{{8.8}{38}{homotopy pushouts}{thm.8.8}{}}
\newlabel{reedy:eq:pushout.compare}{{8.9}{39}{homotopy pushouts}{equation.8.9}{}}
\newlabel{reedy:ex:stupid-simplicial}{{8.10}{39}{}{thm.8.10}{}}
\newlabel{reedy:sec:connected-weights}{{9}{39}{Connected weights}{section.9}{}}
\newlabel{reedy:prop:2/3-SM7}{{9.1}{39}{}{thm.9.1}{}}
\newlabel{reedy:obs:fibrant.constants}{{9.2}{40}{}{thm.9.2}{}}
\newlabel{reedy:cor:connected.weights}{{9.4}{40}{}{thm.9.4}{}}
\newlabel{reedy:sec:simplicial}{{10}{41}{Simplicial model categories and geometric realization}{section.10}{}}
\newlabel{reedy:thm:simp.model.cat}{{10.3}{42}{}{thm.10.3}{}}
\newlabel{reedy:eq:simpmodelreedy}{{10.4}{42}{Simplicial model categories and geometric realization}{equation.10.4}{}}
\newlabel{reedy:obs:geo-filt}{{10.7}{43}{skeletal filtration of geometric realization}{thm.10.7}{}}
\newlabel{reedy:eq:reedy's.po2}{{10.8}{43}{skeletal filtration of geometric realization}{equation.10.8}{}}
\newlabel{reedy:tocindent-1}{0pt}
\newlabel{reedy:tocindent0}{15.01021pt}
\newlabel{reedy:tocindent1}{26.75734pt}
\newlabel{reedy:tocindent2}{0pt}
\newlabel{reedy:tocindent3}{0pt}

%% file: complete_lab.tex
\newlabel{complete:thm:limits-in-limits}{{1.1}{2}{}{thm.1.1}{}}
\newlabel{complete:item:stepI}{{{{(I)}}}{4}{Outline of the proof}{Item.1}{}}
\newlabel{complete:item:stepII}{{{{(II)}}}{4}{Outline of the proof}{Item.2}{}}
\newlabel{complete:sec:weighted}{{2}{5}{Weighted limits in \texorpdfstring {$\protect \qCat _\infty $}{qCat}}{section.2}{}}
\newlabel{complete:sec:terminal}{{3}{6}{Quasi-categories with terminal objects}{section.3}{}}
\newlabel{complete:defn:qCat-X}{{3.1}{6}{}{thm.3.1}{}}
\newlabel{complete:obs:terminal-pres-banal}{{3.2}{6}{}{thm.3.2}{}}
\newlabel{complete:obs:terminals-in-limits-prf}{{3.3}{6}{}{thm.3.3}{}}
\newlabel{complete:rmk:qCat-term-closure-expl}{{3.4}{7}{proof-schema}{thm.3.4}{}}
\newlabel{complete:enum:term-in-lim-p1}{{{{(i)}}}{7}{proof-schema}{Item.13}{}}
\newlabel{complete:enum:term-in-lim-p2}{{{{(ii)}}}{7}{proof-schema}{Item.14}{}}
\newlabel{complete:rec:terminal.char}{{3.5}{7}{terminal objects}{thm.3.5}{}}
\newlabel{complete:lem:rel-lift-terminal-conv}{{3.7}{8}{}{thm.3.7}{}}
\newlabel{complete:lem:rel-lift-terminal-companion}{{3.8}{8}{}{thm.3.8}{}}
\newlabel{complete:lem:qCat_t-closed-products}{{3.9}{8}{}{thm.3.9}{}}
\newlabel{complete:eq:to-solve-prod}{{3.10}{9}{Conical limits in \texorpdfstring {$\protect \qCat _\emptyset $}{qCat}}{equation.3.10}{}}
\newlabel{complete:lem:qCat_t-closed-pullbacks}{{3.11}{9}{}{thm.3.11}{}}
\newlabel{complete:eq:pb-in-qCat}{{3}{9}{Conical limits in \texorpdfstring {$\protect \qCat _\emptyset $}{qCat}}{thm.3.11}{}}
\newlabel{complete:eq:to-solve-pb}{{3.12}{9}{Conical limits in \texorpdfstring {$\protect \qCat _\emptyset $}{qCat}}{equation.3.12}{}}
\newlabel{complete:lem:qCat_t-closed-countable-composite}{{3.13}{10}{}{thm.3.13}{}}
\newlabel{complete:eq:to-solve-tower}{{3.14}{10}{Conical limits in \texorpdfstring {$\protect \qCat _\emptyset $}{qCat}}{equation.3.14}{}}
\newlabel{complete:lem:qCat_t-closed-idempotents}{{3.15}{11}{}{thm.3.15}{}}
\newlabel{complete:thm:terminal-closure}{{3.16}{11}{}{thm.3.16}{}}
\newlabel{complete:sec:absolute}{{4}{12}{Absolute lifting diagrams via terminal objects}{section.4}{}}
\newlabel{complete:eq:obj-qCat-pbshape}{{4.2}{12}{}{equation.4.2}{}}
\newlabel{complete:eq:obj-qCat-pbshape-trans}{{4.3}{12}{}{equation.4.3}{}}
\newlabel{complete:eq:lifting.triangle}{{4.4}{12}{}{equation.4.4}{}}
\newlabel{complete:defn:abs-lift-trans-right-exactness}{{4.5}{13}{}{thm.4.5}{}}
\newlabel{complete:eq:induced-mate}{{4.6}{13}{}{equation.4.6}{}}
\newlabel{complete:obs:comparison-2cell-composite}{{4.8}{13}{}{thm.4.8}{}}
\newlabel{complete:prop:abslifts-in-limits}{{4.9}{14}{}{thm.4.9}{}}
\newlabel{complete:eq:limit-X-pbshape}{{4.10}{14}{Projective cofibrant weighted limits in \texorpdfstring {$\protect \qCat _{r\infty }^{\protect \pbshape }$}{qCat}}{equation.4.10}{}}
\newlabel{complete:eq:limit-X-pbshape-trans}{{4.11}{14}{Projective cofibrant weighted limits in \texorpdfstring {$\protect \qCat _{r\infty }^{\protect \pbshape }$}{qCat}}{equation.4.11}{}}
\newlabel{complete:obs:comma-simp-functor}{{4.12}{15}{}{thm.4.12}{}}
\newlabel{complete:obs:pointwise-exactness}{{4.13}{15}{a pointwise characterisation of right exactness}{thm.4.13}{}}
\newlabel{complete:lem:right-exact-pointwise}{{4.14}{16}{}{thm.4.14}{}}
\newlabel{complete:lem:conical-closure}{{4.15}{16}{}{thm.4.15}{}}
\newlabel{complete:eq:conical-limit-cone}{{4.16}{17}{Projective cofibrant weighted limits in \texorpdfstring {$\protect \qCat _{r\infty }^{\protect \pbshape }$}{qCat}}{equation.4.16}{}}
\newlabel{complete:lem:cotensor-closure}{{4.17}{18}{}{thm.4.17}{}}
\newlabel{complete:eq:cotensor-in-qCat-pb}{{4.18}{18}{Projective cofibrant weighted limits in \texorpdfstring {$\protect \qCat _{r\infty }^{\protect \pbshape }$}{qCat}}{equation.4.18}{}}
\newlabel{complete:eq:induced-into-cotensor}{{4.19}{19}{Projective cofibrant weighted limits in \texorpdfstring {$\protect \qCat _{r\infty }^{\protect \pbshape }$}{qCat}}{equation.4.19}{}}
\newlabel{complete:cor:terminal-cotensor}{{4.20}{20}{}{thm.4.20}{}}
\newlabel{complete:thm:stable}{{4.21}{20}{}{thm.4.21}{}}
\newlabel{complete:sec:algebras}{{5}{22}{Limits and colimits in the quasi-category of algebras}{section.5}{}}
\newlabel{complete:rec:algebras}{{5.1}{22}{homotopy coherent monads and their algebras}{thm.5.1}{}}
\newlabel{complete:obs:algcats-functoriality}{{5.2}{22}{}{thm.5.2}{}}
\newlabel{complete:eq:commutes-with-underlying}{{5.3}{23}{}{equation.5.3}{}}
\newlabel{complete:cor:monadic-colimits}{{5.5}{23}{}{thm.5.5}{}}
\newlabel{complete:cor:monadic-creates-u-split}{{5.6}{23}{}{thm.5.6}{}}
\newlabel{complete:thm:monadic-completeness}{{5.7}{23}{}{thm.5.7}{}}
\newlabel{complete:lem:W-.decomp}{{5.10}{25}{}{thm.5.10}{}}
\newlabel{complete:eq:projective-cell}{{5.11}{25}{Understanding the monadic forgetful functor}{equation.5.11}{}}
\newlabel{complete:cor:W-.decomp}{{5.12}{25}{}{thm.5.12}{}}
\newlabel{complete:eq:tower}{{5.13}{26}{Understanding the monadic forgetful functor}{equation.5.13}{}}
\newlabel{complete:eq:tower.step.pb}{{5.14}{26}{Understanding the monadic forgetful functor}{equation.5.14}{}}
\newlabel{complete:obs:W-.decomp.key}{{5.15}{26}{The key fact}{thm.5.15}{}}
\newlabel{complete:eq:terminal-fact}{{5.16}{26}{The key fact}{equation.5.16}{}}
\newlabel{complete:prop:monadic-terminal}{{5.17}{27}{}{thm.5.17}{}}
\newlabel{complete:eq:monadic-terminal}{{5.18}{27}{}{equation.5.18}{}}
\newlabel{complete:cor:monadic-terminal}{{5.20}{27}{}{thm.5.20}{}}
\newlabel{complete:eq:ind.step.filler.diag}{{5.21}{28}{Creation of terminal objects}{equation.5.21}{}}
\newlabel{complete:eq:lift1}{{5.22}{29}{Creation of limits}{equation.5.22}{}}
\newlabel{complete:thm:monadic-limit}{{5.24}{29}{}{thm.5.24}{}}
\newlabel{complete:eq:monadic-hypothesis}{{5.25}{29}{}{equation.5.25}{}}
\newlabel{complete:eq:monadic-conclusion}{{5.26}{29}{}{equation.5.26}{}}
\newlabel{complete:eq:pres-by-underlying}{{5.27}{30}{}{equation.5.27}{}}
\newlabel{complete:prop:monadic-limit-2}{{5.29}{30}{}{thm.5.29}{}}
\newlabel{complete:obs:expl-alg-desc}{{5.30}{30}{}{thm.5.30}{}}
\newlabel{complete:eq:expl-alg-desc}{{5.31}{30}{}{equation.5.31}{}}
\newlabel{complete:obs:W-.monad}{{5.32}{30}{}{thm.5.32}{}}
\newlabel{complete:eq:pre-commas}{{5.34}{31}{}{equation.5.34}{}}
\newlabel{complete:eq:squares-of-power}{{5.35}{32}{}{equation.5.35}{}}
\newlabel{complete:lem:technical-lem}{{5.36}{32}{}{thm.5.36}{}}
\newlabel{complete:tocindent-1}{0pt}
\newlabel{complete:tocindent0}{15.01021pt}
\newlabel{complete:tocindent1}{20.88377pt}
\newlabel{complete:tocindent2}{0pt}
\newlabel{complete:tocindent3}{0pt}

%% file: equipment_lab.tex
\newlabel{equipment:sec:background}{{2}{6}{Background}{section.2}{}}
\newlabel{equipment:ssec:cosmoi}{{2.1}{6}{\texorpdfstring {$\infty $}{infinity}-cosmoi and their homotopy 2-categories}{subsection.2.1}{}}
\newlabel{equipment:qcat.ctxt.cof.def}{{2.1.1}{6}{$\infty $-cosmos}{thm.2.1.1}{}}
\newlabel{equipment:qcat.ctxt.cof:a}{{{{(a)}}}{6}{$\infty $-cosmos}{Item.1}{}}
\newlabel{equipment:qcat.ctxt.cof:b}{{{{(b)}}}{6}{$\infty $-cosmos}{Item.2}{}}
\newlabel{equipment:qcat.ctxt.cof:c}{{{{(c)}}}{7}{$\infty $-cosmos}{Item.3}{}}
\newlabel{equipment:lem:triv.fib.stab}{{2.1.3}{7}{stability properties of trivial fibrations}{thm.2.1.3}{}}
\newlabel{equipment:itm:triv.fib.stab:a}{{{{(a)}}}{7}{stability properties of trivial fibrations}{Item.4}{}}
\newlabel{equipment:itm:triv.fib.stab:b}{{{{(b)}}}{7}{stability properties of trivial fibrations}{Item.5}{}}
\newlabel{equipment:itm:triv.fib.stab:c}{{{{(c)}}}{7}{stability properties of trivial fibrations}{Item.6}{}}
\newlabel{equipment:defn:closed-cosmoi}{{2.1.5}{7}{cartesian closed $\infty $-cosmoi}{thm.2.1.5}{}}
\newlabel{equipment:defn:sliced-cosmoi}{{2.1.6}{8}{sliced $\infty $-cosmoi}{thm.2.1.6}{}}
\newlabel{equipment:eq:sliced-map}{{2.1.7}{8}{sliced $\infty $-cosmoi}{equation.2.1.7}{}}
\newlabel{equipment:defn:qcat-ctxt-functor}{{2.1.8}{8}{}{thm.2.1.8}{}}
\newlabel{equipment:prop:cartesian-closure}{{2.1.11}{9}{}{thm.2.1.11}{}}
\newlabel{equipment:defn:groupoidal-object}{{2.1.12}{10}{}{thm.2.1.12}{}}
\newlabel{equipment:ssec:sliced}{{2.2}{10}{Sliced homotopy 2-categories}{subsection.2.2}{}}
\newlabel{equipment:defn:two-slices}{{2.2.1}{10}{}{thm.2.2.1}{}}
\newlabel{equipment:prop:smothering-slice-comparison}{{2.2.2}{11}{}{thm.2.2.2}{}}
\newlabel{equipment:cor:slice-2-cat-comparison}{{2.2.3}{11}{}{thm.2.2.3}{}}
\newlabel{equipment:itm:slice-2-cat-comparison-equiv}{{{{(ii)}}}{11}{}{Item.10}{}}
\newlabel{equipment:itm:slice-2-cat-comparison-iso-class}{{{{(iii)}}}{11}{}{Item.11}{}}
\newlabel{equipment:itm:slice-2-cat-comparison-groupoidal}{{{{(iv)}}}{11}{}{Item.12}{}}
\newlabel{equipment:itm:slice-2-cat-comparison-adj}{{{{(v)}}}{11}{}{Item.13}{}}
\newlabel{equipment:rmk:adjunctions-pullback}{{2.2.6}{12}{}{thm.2.2.6}{}}
\newlabel{equipment:ssec:limits}{{2.3}{12}{Simplicial limits modeling comma \texorpdfstring {$\infty $}{infinity}-categories}{subsection.2.3}{}}
\newlabel{equipment:rec:comma}{{2.3.1}{12}{comma $\infty $-categories}{thm.2.3.1}{}}
\newlabel{equipment:eq:comma-as-simp-pullback}{{2.3.2}{13}{comma $\infty $-categories}{equation.2.3.2}{}}
\newlabel{equipment:eq:comma-square}{{2.3.3}{13}{comma $\infty $-categories}{equation.2.3.3}{}}
\newlabel{equipment:rec:ess.unique.1-cell.ind}{{2.3.4}{13}{}{thm.2.3.4}{}}
\newlabel{equipment:lem:equiv-invariance-of-commas}{{2.3.6}{14}{}{thm.2.3.6}{}}
\newlabel{equipment:lem:commas-pullback}{{2.3.7}{14}{}{thm.2.3.7}{}}
\newlabel{equipment:lem:equiv-to-comma}{{2.3.8}{15}{}{thm.2.3.8}{}}
\newlabel{equipment:ex:cone-comma-squares}{{2.3.9}{15}{}{thm.2.3.9}{}}
\newlabel{equipment:prop:induced-2-functor}{{2.3.10}{16}{}{thm.2.3.10}{}}
\newlabel{equipment:ssec:cartesian}{{2.4}{16}{Cartesian fibrations and groupoidal cartesian fibrations}{subsection.2.4}{}}
\newlabel{equipment:itm:weak.cart.i}{{{{(i)}}}{16}{cartesian 2-cells}{Item.17}{}}
\newlabel{equipment:itm:weak.cart.ii}{{{{(ii)}}}{16}{cartesian 2-cells}{Item.18}{}}
\newlabel{equipment:eq:2-cell.to.lift}{{2.4.3}{16}{cartesian fibration}{equation.2.4.3}{}}
\newlabel{equipment:itm:cart.fib.chars.i}{{{{(i)}}}{17}{}{Item.21}{}}
\newlabel{equipment:itm:cart.fib.chars.ii}{{{{(ii)}}}{17}{}{Item.22}{}}
\newlabel{equipment:itm:cart.fib.chars.iii}{{{{(iii)}}}{17}{}{Item.23}{}}
\newlabel{equipment:thm:gpd.cart.fib.chars}{{2.4.5}{17}{}{thm.2.4.5}{}}
\newlabel{equipment:itm:gpd.cart.fib.chars.i}{{{{(i)}}}{17}{}{Item.24}{}}
\newlabel{equipment:itm:gpd.cart.fib.chars.ii}{{{{(ii)}}}{17}{}{Item.25}{}}
\newlabel{equipment:itm:gpd.cart.fib.chars.iii}{{{{(iii)}}}{17}{}{Item.26}{}}
\newlabel{equipment:defn:cartesian-functor}{{2.4.7}{18}{}{thm.2.4.7}{}}
\newlabel{equipment:prop:cart-fib-pullback}{{2.4.8}{18}{\refIV {prop:cart-fib-pullback} and \refIV {cor:groupoidal-pullback}}{thm.2.4.8}{}}
\newlabel{equipment:sec:modules}{{3}{19}{Modules between \texorpdfstring {$\infty $}{infinity}-categories}{section.3}{}}
\newlabel{equipment:itm:cartesian-on-the-right}{{{{(i)}}}{19}{}{Item.29}{}}
\newlabel{equipment:itm:cocartesian-on-the-left}{{{{(ii)}}}{19}{}{Item.30}{}}
\newlabel{equipment:itm:groupoidal-fibers}{{{{(iii)}}}{19}{}{Item.31}{}}
\newlabel{equipment:lem:disc-cart-on-right}{{3.1.3}{20}{}{thm.3.1.3}{}}
\newlabel{equipment:prop:hom-is-a-module}{{3.1.4}{21}{}{thm.3.1.4}{}}
\newlabel{equipment:lem:sliced-adjunction-on-the-right}{{3.1.5}{21}{}{thm.3.1.5}{}}
\newlabel{equipment:itm:sliced-right-defn}{{{{(i)}}}{21}{}{Item.32}{}}
\newlabel{equipment:itm:sliced-right-fibred-adj}{{{{(ii)}}}{21}{}{Item.33}{}}
\newlabel{equipment:itm:sliced-right-adj}{{{{(iii)}}}{21}{}{Item.34}{}}
\newlabel{equipment:prop:two-sided-pullback}{{3.1.6}{22}{}{thm.3.1.6}{}}
\newlabel{equipment:cor:comma-is-a-module}{{3.1.7}{23}{}{thm.3.1.7}{}}
\newlabel{equipment:lem:two-sided-groupoidal-cells}{{3.1.9}{23}{}{thm.3.1.9}{}}
\newlabel{equipment:defn:horiz-comp}{{3.1.10}{24}{horizontal composition of isofibrations over products}{thm.3.1.10}{}}
\newlabel{equipment:lem:module-pullback-cartesian}{{3.1.12}{25}{}{thm.3.1.12}{}}
\newlabel{equipment:ex:modules-do-not-compose}{{3.1.13}{25}{}{thm.3.1.13}{}}
\newlabel{equipment:lem:span-maps-cartesian}{{3.2.1}{26}{}{thm.3.2.1}{}}
\newlabel{equipment:eq:span-map}{{3.2.2}{26}{}{equation.3.2.2}{}}
\newlabel{equipment:prop:two-sided-yoneda}{{3.2.4}{26}{}{thm.3.2.4}{}}
\newlabel{equipment:defn:1-cat-of-modules}{{3.2.5}{27}{}{thm.3.2.5}{}}
\newlabel{equipment:lem:full-yoneda}{{3.2.6}{28}{}{thm.3.2.6}{}}
\newlabel{equipment:defn:module-map}{{3.2.7}{28}{}{thm.3.2.7}{}}
\newlabel{equipment:lem:module-map-pullback-equivalence}{{3.2.8}{29}{}{thm.3.2.8}{}}
\newlabel{equipment:eq:pullback-for-UP}{{3.2.9}{29}{}{equation.3.2.9}{}}
\newlabel{equipment:lem:equiv-of-modules}{{3.3.1}{29}{}{thm.3.3.1}{}}
\newlabel{equipment:lem:rep-mod-char}{{3.3.2}{30}{}{thm.3.3.2}{}}
\newlabel{equipment:eq:r-equiv-defn}{{3.3.3}{30}{Equivalence of modules}{equation.3.3.3}{}}
\newlabel{equipment:eq:e-equiv-defn}{{3.3.4}{30}{Equivalence of modules}{equation.3.3.4}{}}
\newlabel{equipment:lem:qcat-pointwise-rep}{{3.3.5}{31}{}{thm.3.3.5}{}}
\newlabel{equipment:eq:lift-term-in-fibre}{{3.3.6}{31}{Equivalence of modules}{equation.3.3.6}{}}
\newlabel{equipment:cor:fibrewise-rep}{{3.3.7}{32}{}{thm.3.3.7}{}}
\newlabel{equipment:ex:adjoint-equivalence-mods}{{3.3.8}{32}{}{thm.3.3.8}{}}
\newlabel{equipment:sec:virtual}{{4}{33}{The virtual equipment of modules}{section.4}{}}
\newlabel{equipment:defn:double-iso-spans}{{4.1.1}{34}{the double category of isofibrations}{thm.4.1.1}{}}
\newlabel{equipment:eq:generic-cell}{{4.1.3}{35}{virtual double category}{equation.4.1.3}{}}
\newlabel{equipment:prop:modules-virtual}{{4.1.4}{35}{}{thm.4.1.4}{}}
\newlabel{equipment:defn:composable-modules}{{4.1.5}{36}{composable modules}{thm.4.1.5}{}}
\newlabel{equipment:obs:cells-above-a-comma}{{4.1.6}{36}{}{thm.4.1.6}{}}
\newlabel{equipment:eq:simplicial-pullback-of-module}{{4.2.1}{37}{The virtual equipment of modules}{equation.4.2.1}{}}
\newlabel{equipment:prop:cartesian-pullback-cells}{{4.2.2}{37}{}{thm.4.2.2}{}}
\newlabel{equipment:eq:cartesian-cell-UP}{{4.2.3}{37}{}{equation.4.2.3}{}}
\newlabel{equipment:prop:arrows-are-units}{{4.2.4}{38}{}{thm.4.2.4}{}}
\newlabel{equipment:defn:virtual-equipment}{{4.2.5}{39}{\cite [\S 7]{CruttwellShulman:2010au}}{thm.4.2.5}{}}
\newlabel{equipment:thm:virtual-equipment}{{4.2.6}{39}{}{thm.4.2.6}{}}
\newlabel{equipment:eq:multi-domain}{{4.3.3}{40}{composition of modules}{equation.4.3.3}{}}
\newlabel{equipment:eq:generic-composite}{{4.3.4}{40}{composition of modules}{equation.4.3.4}{}}
\newlabel{equipment:obs:unary-comp}{{4.3.5}{40}{nullary and unary composites}{thm.4.3.5}{}}
\newlabel{equipment:obs:comp-assoc}{{4.3.6}{41}{associativity of composition}{thm.4.3.6}{}}
\newlabel{equipment:obs:comp-cancel}{{4.3.7}{41}{left cancelation of composites}{thm.4.3.7}{}}
\newlabel{equipment:obs:comp-proof-strategy}{{4.3.8}{41}{}{thm.4.3.8}{}}
\newlabel{equipment:lem:comp-with-unit}{{4.3.9}{41}{composites with units}{thm.4.3.9}{}}
\newlabel{equipment:eq:comp-with-unit}{{4.3.10}{42}{composites with units}{equation.4.3.10}{}}
\newlabel{equipment:obs:comp-with-unit-over-comma}{{4.3.11}{42}{}{thm.4.3.11}{}}
\newlabel{equipment:defn:unit-cells}{{4.3.12}{42}{unit cells}{thm.4.3.12}{}}
\newlabel{equipment:lem:comp-with-unit-cells}{{4.3.13}{43}{composite with unit cells}{thm.4.3.13}{}}
\newlabel{equipment:defn:companion-conjoint}{{4.4.1}{44}{}{thm.4.4.1}{}}
\newlabel{equipment:thm:companion-conjoint}{{4.4.2}{45}{}{thm.4.4.2}{}}
\newlabel{equipment:cor:companion-conjoint}{{4.4.3}{46}{}{thm.4.4.3}{}}
\newlabel{equipment:lem:fibred-adjoint-reflection}{{4.4.4}{47}{}{thm.4.4.4}{}}
\newlabel{equipment:lem:one-sided-rep-comp}{{4.4.5}{47}{}{thm.4.4.5}{}}
\newlabel{equipment:cor:two-sided-rep-comp}{{4.4.7}{48}{}{thm.4.4.7}{}}
\newlabel{equipment:ex:comp-of-rep}{{4.4.8}{48}{}{thm.4.4.8}{}}
\newlabel{equipment:ex:two-sided-comma-as-composite}{{4.4.9}{49}{}{thm.4.4.9}{}}
\newlabel{equipment:lem:module-as-rep-composite}{{4.4.10}{49}{}{thm.4.4.10}{}}
\newlabel{equipment:lem:yoneda-embedding}{{4.4.11}{50}{}{thm.4.4.11}{}}
\newlabel{equipment:rmk:yoneda-embedding}{{4.4.12}{50}{}{thm.4.4.12}{}}
\newlabel{equipment:sec:pointwise-kan}{{5}{50}{Pointwise Kan extensions}{section.5}{}}
\newlabel{equipment:ssec:exact-squares}{{5.1}{51}{Exact squares}{subsection.5.1}{}}
\newlabel{equipment:defn:exact-square}{{5.1.1}{51}{exact squares}{thm.5.1.1}{}}
\newlabel{equipment:lem:exact-squares-comp}{{5.1.2}{52}{composites of exact squares}{thm.5.1.2}{}}
\newlabel{equipment:lem:comma-exact}{{5.1.3}{53}{}{thm.5.1.3}{}}
\newlabel{equipment:lem:pullback-exact}{{5.1.4}{53}{}{thm.5.1.4}{}}
\newlabel{equipment:lem:product-exact-square}{{5.1.5}{54}{}{thm.5.1.5}{}}
\newlabel{equipment:lem:product-comma-exact}{{5.1.6}{55}{}{thm.5.1.6}{}}
\newlabel{equipment:ssec:pointwise-kan}{{5.2}{56}{Pointwise Kan Extensions}{subsection.5.2}{}}
\newlabel{equipment:defn:right-extension-mod}{{5.2.1}{56}{right extension of modules}{thm.5.2.1}{}}
\newlabel{equipment:lem:mod-extensions-in-2-cat}{{5.2.2}{56}{}{thm.5.2.2}{}}
\newlabel{equipment:prop:pointwise-kan}{{5.2.4}{57}{}{thm.5.2.4}{}}
\newlabel{equipment:eq:pointwise-right-extension}{{5.2.5}{57}{}{equation.5.2.5}{}}
\newlabel{equipment:itm:pointwise-exact-stable}{{{{(i)}}}{57}{}{Item.40}{}}
\newlabel{equipment:itm:pointwise-comma-stable}{{{{(ii)}}}{57}{}{Item.41}{}}
\newlabel{equipment:itm:pointwise-in-mod}{{{{(iii)}}}{57}{}{Item.42}{}}
\newlabel{equipment:itm:pointwise-exact-stable-in-mod}{{{{(iv)}}}{57}{}{Item.43}{}}
\newlabel{equipment:defn:fully-faithful}{{5.2.7}{59}{fully faithful}{thm.5.2.7}{}}
\newlabel{equipment:ssec:pointwise-kan-closed}{{5.3}{60}{Pointwise Kan extensions in a cartesian closed \texorpdfstring {$\infty $}{infinity}-cosmos}{subsection.5.3}{}}
\newlabel{equipment:prop:ran-abs-lifting}{{5.3.1}{61}{}{thm.5.3.1}{}}
\newlabel{equipment:eq:pointwise-kan-abs-lifting}{{5.3.2}{61}{}{equation.5.3.2}{}}
\newlabel{equipment:eq:limit-abs-lifting}{{5.3.3}{62}{Pointwise Kan extensions in a cartesian closed \texorpdfstring {$\infty $}{infinity}-cosmos}{equation.5.3.3}{}}
\newlabel{equipment:prop:limits-as-pointwise-kan}{{5.3.4}{62}{}{thm.5.3.4}{}}
\newlabel{equipment:defn:initial-final}{{5.3.5}{62}{initial/final functor}{thm.5.3.5}{}}
\newlabel{equipment:prop:beck-chevalley}{{5.3.9}{64}{Beck-Chevalley condition}{thm.5.3.9}{}}
\newlabel{equipment:rmk:qcat-derivator}{{5.3.10}{65}{on derivators for (co)complete quasi-categories}{thm.5.3.10}{}}
\newlabel{equipment:eq:qcat-derivator}{{5.3.11}{66}{on derivators for (co)complete quasi-categories}{equation.5.3.11}{}}
\newlabel{equipment:eq:mod-right-ext}{{5.3.12}{67}{on derivators for (co)complete quasi-categories}{equation.5.3.12}{}}
\newlabel{equipment:tocindent-1}{0pt}
\newlabel{equipment:tocindent0}{15.01021pt}
\newlabel{equipment:tocindent1}{20.88377pt}
\newlabel{equipment:tocindent2}{34.49158pt}
\newlabel{equipment:tocindent3}{0pt}

%% file: yoneda_lab.tex
\newlabel{yoneda:sec:cosmoi}{{2}{5}{\texorpdfstring {$\infty $}{infinity}-cosmoi}{section.2}{}}
\newlabel{yoneda:ssec:cosmoi}{{2.1}{6}{\texorpdfstring {$\infty $}{infinity}-cosmoi and functors}{subsection.2.1}{}}
\newlabel{yoneda:qcat.ctxt.def}{{2.1.1}{6}{$\infty $-cosmos}{thm.2.1.1}{}}
\newlabel{yoneda:qcat.ctxt:a}{{{{(a)}}}{6}{$\infty $-cosmos}{Item.4}{}}
\newlabel{yoneda:qcat.ctxt:b}{{{{(b)}}}{6}{$\infty $-cosmos}{Item.5}{}}
\newlabel{yoneda:qcat.ctxt:d}{{{{(c)}}}{6}{$\infty $-cosmos}{Item.6}{}}
\newlabel{yoneda:qcat.ctxt:f}{{{{(d)}}}{6}{$\infty $-cosmos}{Item.7}{}}
\newlabel{yoneda:qcat.ctxt:e}{{{{(e)}}}{6}{$\infty $-cosmos}{Item.8}{}}
\newlabel{yoneda:qcat.ctxt.products}{{2.1.2}{6}{}{thm.2.1.2}{}}
\newlabel{yoneda:qcat.ctxt.clarif}{{2.1.3}{6}{}{thm.2.1.3}{}}
\newlabel{yoneda:ex:qcat-qcat-ctxt}{{2.1.4}{7}{the $\infty $-cosmos of quasi-categories}{thm.2.1.4}{}}
\newlabel{yoneda:rec:equiv-and-K}{{2.1.5}{7}{}{thm.2.1.5}{}}
\newlabel{yoneda:rec:equiv-and-K.a}{{{{(a)}}}{7}{}{Item.9}{}}
\newlabel{yoneda:rec:equiv-and-K.b}{{{{(b)}}}{7}{}{Item.10}{}}
\newlabel{yoneda:rec:equiv-and-K.c}{{{{(c)}}}{7}{}{Item.11}{}}
\newlabel{yoneda:rec:equiv-and-K.d}{{{{(d)}}}{7}{}{Item.12}{}}
\newlabel{yoneda:lem:Brown.fact}{{2.1.6}{7}{Brown's factorisation lemma}{thm.2.1.6}{}}
\newlabel{yoneda:obs:Brown.fact.cofibrant}{{2.1.7}{8}{}{thm.2.1.7}{}}
\newlabel{yoneda:lem:isofib-rep}{{2.1.8}{8}{}{thm.2.1.8}{}}
\newlabel{yoneda:defn:qcat-ctxt-functor}{{2.1.9}{9}{}{thm.2.1.9}{}}
\newlabel{yoneda:prop:representable-functors}{{2.1.10}{9}{}{thm.2.1.10}{}}
\newlabel{yoneda:ex:sliced.contexts}{{2.1.11}{9}{sliced $\infty $-cosmoi}{thm.2.1.11}{}}
\newlabel{yoneda:eq:sliced-map}{{2.1.12}{9}{sliced $\infty $-cosmoi}{equation.2.1.12}{}}
\newlabel{yoneda:prop:pullback-functors}{{2.1.13}{10}{pulling back between sliced $\infty $-cosmoi}{thm.2.1.13}{}}
\newlabel{yoneda:ssec:cosmoi-examples}{{2.2}{11}{Examples of $\infty $-cosmoi}{subsection.2.2}{}}
\newlabel{yoneda:lem:model-categories-cosmoi}{{2.2.1}{11}{}{thm.2.2.1}{}}
\newlabel{yoneda:obs:duals-of-cosmoi}{{2.2.2}{11}{}{thm.2.2.2}{}}
\newlabel{yoneda:prop:model-categories-cosmoi}{{2.2.3}{11}{}{thm.2.2.3}{}}
\newlabel{yoneda:ex:cat-cosmos}{{2.2.4}{11}{the $\infty $-cosmos of categories}{thm.2.2.4}{}}
\newlabel{yoneda:ex:CSS-cosmos}{{2.2.5}{12}{the $\infty $-cosmos of complete Segal spaces}{thm.2.2.5}{}}
\newlabel{yoneda:ex:other-CSS-functor}{{2.2.6}{12}{}{thm.2.2.6}{}}
\newlabel{yoneda:ex:segal-cosmos}{{2.2.7}{13}{the $\infty $-cosmos of Segal categories}{thm.2.2.7}{}}
\newlabel{yoneda:ex:marked-cosmos}{{2.2.8}{14}{the $\infty $-cosmos of marked simplicial sets}{thm.2.2.8}{}}
\newlabel{yoneda:prop:rezk-cosmos}{{2.2.9}{14}{the $\infty $-cosmos of Rezk objects}{thm.2.2.9}{}}
\newlabel{yoneda:ex:theta-n-cosmos}{{2.2.10}{16}{the $\infty $-cosmos of $\theta _n$-spaces}{thm.2.2.10}{}}
\newlabel{yoneda:eq:CSS-theta-n-adj}{{2.2.11}{16}{the $\infty $-cosmos of $\theta _n$-spaces}{equation.2.2.11}{}}
\newlabel{yoneda:sec:abstract}{{3}{17}{2-category theory in an \texorpdfstring {$\infty $}{infinity}-cosmos}{section.3}{}}
\newlabel{yoneda:ssec:htpy-2-cat}{{3.1}{17}{The homotopy 2-category of an \texorpdfstring {$\infty $}{infinity}-cosmos}{subsection.3.1}{}}
\newlabel{yoneda:rec:trivial-fibration}{{3.1.3}{18}{equivalences, isofibrations, and surjective equivalences}{thm.3.1.3}{}}
\newlabel{yoneda:lem:isofib.are.representably.so}{{3.1.4}{18}{}{thm.3.1.4}{}}
\newlabel{yoneda:cor:isofib.are.representably.so}{{3.1.5}{19}{}{thm.3.1.5}{}}
\newlabel{yoneda:obs:representing-isos}{{3.1.6}{19}{functors representing isomorphic 2-cells}{thm.3.1.6}{}}
\newlabel{yoneda:lem:iso-to-weak-equiv}{{3.1.7}{19}{}{thm.3.1.7}{}}
\newlabel{yoneda:prop:equiv.are.weak.equiv}{{3.1.8}{19}{}{thm.3.1.8}{}}
\newlabel{yoneda:obs:bi-coreflection}{{3.1.9}{20}{cofibrant replacement as a {\em bi-coreflection\/}}{thm.3.1.9}{}}
\newlabel{yoneda:ssec:abstract}{{3.2}{20}{Abstract homotopy 2-categories}{subsection.3.2}{}}
\newlabel{yoneda:defn:comma-object}{{3.2.1}{21}{comma object}{thm.3.2.1}{}}
\newlabel{yoneda:eq:comma-square}{{3.2.2}{21}{comma object}{equation.3.2.2}{}}
\newlabel{yoneda:defn:iso-comma-object}{{3.2.3}{22}{iso-comma object}{thm.3.2.3}{}}
\newlabel{yoneda:eq:iso-comma}{{3.2.4}{22}{iso-comma object}{equation.3.2.4}{}}
\newlabel{yoneda:defn:abstract-htpy-2-cat}{{3.2.5}{22}{abstract homotopy 2-category}{thm.3.2.5}{}}
\newlabel{yoneda:ssec:htpy-2-cat-is-abstract}{{3.3}{23}{Comma and iso-comma objects in a homotopy 2-category}{subsection.3.3}{}}
\newlabel{yoneda:lem:htpy-2-cat-commas}{{3.3.1}{23}{}{thm.3.3.1}{}}
\newlabel{yoneda:eq:comma-as-simp-pullback}{{3.3.2}{23}{}{equation.3.3.2}{}}
\newlabel{yoneda:lem:htpy-2-cat-iso-commas}{{3.3.3}{23}{}{thm.3.3.3}{}}
\newlabel{yoneda:eq:iso-comma-as-simp-pullback}{{3.3.4}{24}{}{equation.3.3.4}{}}
\newlabel{yoneda:obs:legs-are-isofibrations}{{3.3.6}{25}{}{thm.3.3.6}{}}
\newlabel{yoneda:ssec:comma}{{3.4}{25}{Stability and uniqueness of comma and iso-comma constructions}{subsection.3.4}{}}
\newlabel{yoneda:defn:span-isofib}{{3.4.1}{25}{}{thm.3.4.1}{}}
\newlabel{yoneda:itm:isofib-right}{{{{(ii)}}}{25}{}{Item.22}{}}
\newlabel{yoneda:itm:isofib-left}{{{{(iii)}}}{25}{}{Item.23}{}}
\newlabel{yoneda:lem:comma-span-isofib}{{3.4.2}{25}{}{thm.3.4.2}{}}
\newlabel{yoneda:eq:span-map}{{3.4.4}{26}{}{equation.3.4.4}{}}
\newlabel{yoneda:lem:equivalence-of-spans}{{3.4.5}{26}{}{thm.3.4.5}{}}
\newlabel{yoneda:obs:ess.unique.1-cell.ind}{{3.4.6}{27}{essential uniqueness of induced 1-cells}{thm.3.4.6}{}}
\newlabel{yoneda:rec:comma-uniqueness}{{3.4.7}{27}{}{thm.3.4.7}{}}
\newlabel{yoneda:lem:comma-uniqueness}{{3.4.8}{28}{stability of (iso-)comma objects under equivalence}{thm.3.4.8}{}}
\newlabel{yoneda:cor:comma-uniqueness}{{3.4.9}{29}{}{thm.3.4.9}{}}
\newlabel{yoneda:lem:comma-composition}{{3.4.10}{29}{}{thm.3.4.10}{}}
\newlabel{yoneda:lem:comma-cancelation}{{3.4.11}{30}{}{thm.3.4.11}{}}
\newlabel{yoneda:lem:restricted-commas-pullback}{{3.4.12}{31}{}{thm.3.4.12}{}}
\newlabel{yoneda:ssec:pullback}{{3.5}{32}{Iso-commas and pullbacks}{subsection.3.5}{}}
\newlabel{yoneda:obs:slice-2-cat}{{3.5.1}{32}{slice 2-categories}{thm.3.5.1}{}}
\newlabel{yoneda:eq:1-cell.in.slice}{{3.5.2}{32}{slice 2-categories}{equation.3.5.2}{}}
\newlabel{yoneda:lem:fibred-equivalence}{{3.5.3}{32}{}{thm.3.5.3}{}}
\newlabel{yoneda:defn:pullback}{{3.5.4}{33}{pullback}{thm.3.5.4}{}}
\newlabel{yoneda:rmk:generalized-1-cell-induction}{{3.5.5}{33}{}{thm.3.5.5}{}}
\newlabel{yoneda:lem:pullbacks-exist}{{3.5.6}{34}{}{thm.3.5.6}{}}
\newlabel{yoneda:defn:pullback-of-an-isofib}{{3.5.7}{34}{}{thm.3.5.7}{}}
\newlabel{yoneda:ex:pullbacks-from-iso-commas}{{3.5.8}{34}{}{thm.3.5.8}{}}
\newlabel{yoneda:ex:old-pullbacks}{{3.5.9}{35}{}{thm.3.5.9}{}}
\newlabel{yoneda:lem:pullback-comp}{{3.5.10}{35}{composition of pullbacks}{thm.3.5.10}{}}
\newlabel{yoneda:ssec:smothering}{{3.6}{36}{Smothering 2-functors and adjunctions}{subsection.3.6}{}}
\newlabel{yoneda:eq:pbk-1-cell}{{3.6.3}{36}{}{equation.3.6.3}{}}
\newlabel{yoneda:eq:pbk-2-cell}{{3.6.4}{37}{}{equation.3.6.4}{}}
\newlabel{yoneda:lem:pullback-smothering}{{3.6.5}{37}{}{thm.3.6.5}{}}
\newlabel{yoneda:cor:adjunctions-pullback}{{3.6.6}{37}{}{thm.3.6.6}{}}
\newlabel{yoneda:cor:pullback-uniqueness}{{3.6.7}{37}{}{thm.3.6.7}{}}
\newlabel{yoneda:defn:iso-comma-slice-2-cat}{{3.6.8}{38}{}{thm.3.6.8}{}}
\newlabel{yoneda:eq:iso-comma-1-cell}{{3.6.9}{38}{}{equation.3.6.9}{}}
\newlabel{yoneda:eq:iso-comma-2-cell}{{3.6.10}{38}{}{equation.3.6.10}{}}
\newlabel{yoneda:lem:iso-comma-smothering}{{3.6.11}{38}{}{thm.3.6.11}{}}
\newlabel{yoneda:cor:adjunctions-iso-comma}{{3.6.12}{39}{}{thm.3.6.12}{}}
\newlabel{yoneda:obs:adj.rep.def}{{3.6.13}{39}{adjunctions are representably defined}{thm.3.6.13}{}}
\newlabel{yoneda:obs:RARI}{{3.6.14}{39}{right adjoint right inverse}{thm.3.6.14}{}}
\newlabel{yoneda:obs:gen.2-cat.results}{{3.6.15}{40}{}{thm.3.6.15}{}}
\newlabel{yoneda:sec:cartesian}{{4}{40}{Cartesian fibrations}{section.4}{}}
\newlabel{yoneda:ssec:cartesian}{{4.1}{40}{Cartesian fibrations}{subsection.4.1}{}}
\newlabel{yoneda:defn:weak.cart.2-cell}{{4.1.1}{40}{cartesian 2-cells}{thm.4.1.1}{}}
\newlabel{yoneda:itm:weak.cart.i}{{{{(i)}}}{40}{cartesian 2-cells}{Item.38}{}}
\newlabel{yoneda:itm:weak.cart.ii}{{{{(ii)}}}{40}{cartesian 2-cells}{Item.39}{}}
\newlabel{yoneda:obs:cart.iso.stab}{{4.1.2}{41}{isomorphism stability}{thm.4.1.2}{}}
\newlabel{yoneda:obs:weak.cart.conserv}{{4.1.3}{41}{more conservativity}{thm.4.1.3}{}}
\newlabel{yoneda:defn:cart-fib}{{4.1.4}{41}{cartesian fibration}{thm.4.1.4}{}}
\newlabel{yoneda:eq:2-cell.to.lift}{{4.1.5}{41}{cartesian fibration}{equation.4.1.5}{}}
\newlabel{yoneda:obs:weak.cart.lift.unique}{{4.1.6}{41}{uniqueness of cartesian lifts}{thm.4.1.6}{}}
\newlabel{yoneda:prop:cart-fib-comp}{{4.1.7}{42}{composites of cartesian fibrations}{thm.4.1.7}{}}
\newlabel{yoneda:ntn:for.commas}{{4.1.8}{42}{}{thm.4.1.8}{}}
\newlabel{yoneda:lem:adj.for.commas}{{4.1.9}{43}{}{thm.4.1.9}{}}
\newlabel{yoneda:thm:cart.fib.chars}{{4.1.10}{44}{}{thm.4.1.10}{}}
\newlabel{yoneda:itm:cart.fib.chars.i}{{{{(i)}}}{44}{}{Item.42}{}}
\newlabel{yoneda:itm:cart.fib.chars.ii}{{{{(ii)}}}{44}{}{Item.43}{}}
\newlabel{yoneda:eq:cartesian.fib.adj}{{4.1.11}{44}{}{equation.4.1.11}{}}
\newlabel{yoneda:itm:cart.fib.chars.iii}{{{{(iii)}}}{44}{}{Item.44}{}}
\newlabel{yoneda:eq:cartesian.isosect.adj}{{4.1.12}{44}{}{equation.4.1.12}{}}
\newlabel{yoneda:obs:fibered-unit-iso}{{4.1.13}{44}{}{thm.4.1.13}{}}
\newlabel{yoneda:obs:constructing-weak-lifts}{{4.1.14}{45}{}{thm.4.1.14}{}}
\newlabel{yoneda:cor:cart-fib-rep}{{4.1.15}{45}{}{thm.4.1.15}{}}
\newlabel{yoneda:ex:domain-fibration}{{4.1.16}{45}{}{thm.4.1.16}{}}
\newlabel{yoneda:obs:domain-cart-lifts}{{4.1.17}{46}{cartesian lifts for domain projections}{thm.4.1.17}{}}
\newlabel{yoneda:ex:general-domain-projection}{{4.1.18}{47}{}{thm.4.1.18}{}}
\newlabel{yoneda:prop:bifibration-adjunction}{{4.1.20}{47}{}{thm.4.1.20}{}}
\newlabel{yoneda:eq:pseudo-slice-adj}{{4.1.21}{48}{Cartesian fibrations}{equation.4.1.21}{}}
\newlabel{yoneda:eq:RARI-lifting}{{4.1.24}{49}{right adjoint right inverse as a lifting property}{equation.4.1.24}{}}
\newlabel{yoneda:cor:lurie-cartesian}{{4.1.24}{50}{}{thm.4.1.24}{}}
\newlabel{yoneda:eq:lurie-cartesian}{{4.1.25}{50}{}{equation.4.1.25}{}}
\newlabel{yoneda:ssec:groupoidal}{{4.2}{51}{Groupoidal cartesian fibrations}{subsection.4.2}{}}
\newlabel{yoneda:defn:groupoidal-object}{{4.2.1}{51}{groupoidal objects}{thm.4.2.1}{}}
\newlabel{yoneda:lem:groupoidal.conserv}{{4.2.2}{51}{}{thm.4.2.2}{}}
\newlabel{yoneda:cor:groupoidal-cart-fib-rep}{{4.2.4}{51}{}{thm.4.2.4}{}}
\newlabel{yoneda:prop:groupoidal-ext-char}{{4.2.5}{51}{}{thm.4.2.5}{}}
\newlabel{yoneda:lem:gpd-cart-cancel}{{4.2.6}{52}{}{thm.4.2.6}{}}
\newlabel{yoneda:prop:iso-groupoidal-char}{{4.2.7}{52}{}{thm.4.2.7}{}}
\newlabel{yoneda:obs:biterminal}{{4.2.10}{53}{biterminal objects in the homotopy 2-category}{thm.4.2.10}{}}
\newlabel{yoneda:ex:groupoidal-representable}{{4.2.11}{54}{}{thm.4.2.11}{}}
\newlabel{yoneda:sec:stability}{{5}{54}{Cartesian functors and pullbacks of cartesian fibrations}{section.5}{}}
\newlabel{yoneda:itm:to-do-cart-cells}{{{{(i)}}}{54}{Cartesian functors and pullbacks of cartesian fibrations}{Item.45}{}}
\newlabel{yoneda:itm:to-do-pullback}{{{{(ii)}}}{54}{Cartesian functors and pullbacks of cartesian fibrations}{Item.46}{}}
\newlabel{yoneda:ssec:cart-fun}{{5.1}{54}{Cartesian functors}{subsection.5.1}{}}
\newlabel{yoneda:eq:cartesian-functor-square}{{5.1.2}{54}{}{equation.5.1.2}{}}
\newlabel{yoneda:thm:cart.fun.chars}{{5.1.4}{55}{}{thm.5.1.4}{}}
\newlabel{yoneda:itm:cart.fun.chars.i}{{{{(i)}}}{55}{}{Item.47}{}}
\newlabel{yoneda:itm:cart.fun.chars.ii}{{{{(ii)}}}{55}{}{Item.48}{}}
\newlabel{yoneda:eq:cart-fun-iso.ii}{{5.1.5}{55}{}{equation.5.1.5}{}}
\newlabel{yoneda:itm:cart.fun.chars.iii}{{{{(iii)}}}{55}{}{Item.49}{}}
\newlabel{yoneda:eq:cart-fun-iso.iii}{{5.1.6}{55}{}{equation.5.1.6}{}}
\newlabel{yoneda:cor:radj-cartesian}{{5.1.7}{57}{}{thm.5.1.7}{}}
\newlabel{yoneda:lem:cart-arrows-compose}{{5.1.8}{58}{}{thm.5.1.8}{}}
\newlabel{yoneda:lem:cart-arrows-cancel}{{5.1.9}{58}{}{thm.5.1.9}{}}
\newlabel{yoneda:ssec:pullback-stability}{{5.2}{59}{Pullback stability}{subsection.5.2}{}}
\newlabel{yoneda:prop:cart-fib-pullback}{{5.2.1}{59}{pullbacks of cartesian fibrations}{thm.5.2.1}{}}
\newlabel{yoneda:cor:groupoidal-pullback}{{5.2.2}{61}{}{thm.5.2.2}{}}
\newlabel{yoneda:sec:yoneda}{{6}{62}{The Yoneda lemma}{section.6}{}}
\newlabel{yoneda:thm:yoneda}{{6.0.1}{62}{Yoneda lemma}{thm.6.0.1}{}}
\newlabel{yoneda:eq:evaluation-comma-adjunction}{{6.0.2}{62}{The Yoneda lemma}{equation.6.0.2}{}}
\newlabel{yoneda:ssec:pullback-functoriality}{{6.1}{62}{Functoriality of pullbacks of cartesian fibrations}{subsection.6.1}{}}
\newlabel{yoneda:defn:pbshape-2-cat}{{6.1.1}{63}{}{thm.6.1.1}{}}
\newlabel{yoneda:eq:1-cell-in-pbshape}{{6.1.2}{63}{}{equation.6.1.2}{}}
\newlabel{yoneda:eq:cartesian-weak-pullback}{{6.1.3}{63}{Functoriality of pullbacks of cartesian fibrations}{equation.6.1.3}{}}
\newlabel{yoneda:defn:square-2-cat}{{6.1.4}{63}{}{thm.6.1.4}{}}
\newlabel{yoneda:eq:1-cell-in-square}{{6.1.5}{63}{}{equation.6.1.5}{}}
\newlabel{yoneda:prop:cartesian-smothering}{{6.1.6}{64}{}{thm.6.1.6}{}}
\newlabel{yoneda:ssec:proof-yoneda}{{6.2}{64}{Proof of the Yoneda lemma}{subsection.6.2}{}}
\newlabel{yoneda:eq:1-cell-in-lax-comma}{{6.2.2}{64}{lax slice 2-category}{equation.6.2.2}{}}
\newlabel{yoneda:obs:lax-slice-adjunction}{{6.2.3}{65}{}{thm.6.2.3}{}}
\newlabel{yoneda:eq:lax-slice-left-adjoint}{{6.2.4}{65}{}{equation.6.2.4}{}}
\newlabel{yoneda:obs:pullback-2-category}{{6.2.5}{65}{}{thm.6.2.5}{}}
\newlabel{yoneda:lem:rep-cart-cells}{{6.2.6}{65}{}{thm.6.2.6}{}}
\newlabel{yoneda:eq:strange-right-adjoint}{{6.2.7}{66}{Proof of the Yoneda lemma}{equation.6.2.7}{}}
\newlabel{yoneda:eq:unrestricted-adjunction}{{6.2.8}{66}{Proof of the Yoneda lemma}{equation.6.2.8}{}}
\newlabel{yoneda:lem:codomain-restriction}{{6.2.9}{66}{}{thm.6.2.9}{}}
\newlabel{yoneda:eq:chi-lifting-phi}{{6.2.10}{67}{Proof of the Yoneda lemma}{equation.6.2.10}{}}
\newlabel{yoneda:eq:restricted-adjunction}{{6.2.11}{67}{Proof of the Yoneda lemma}{equation.6.2.11}{}}
\newlabel{yoneda:lem:last-yoneda-step}{{6.2.12}{67}{}{thm.6.2.12}{}}
\newlabel{yoneda:cor:groupoidal-yoneda}{{6.2.13}{68}{}{thm.6.2.13}{}}
\newlabel{yoneda:sec:appendix}{{7}{68}{Appendix: proof of Theorem~\ref {thm:cart.fib.chars}}{section.7}{}}
\newlabel{yoneda:itm:cart.fib.chars.i'}{{{{(i)}}}{68}{}{Item.50}{}}
\newlabel{yoneda:itm:cart.fib.chars.ii'}{{{{(ii)}}}{68}{}{Item.51}{}}
\newlabel{yoneda:eq:cartesian.fib.adj'}{{7.0.1}{69}{}{equation.7.0.1}{}}
\newlabel{yoneda:itm:cart.fib.chars.iii'}{{{{(iii)}}}{69}{}{Item.52}{}}
\newlabel{yoneda:eq:cartesian.isosect.adj'}{{7.0.2}{69}{}{equation.7.0.2}{}}
\newlabel{yoneda:ntn:adj.for.commas}{{7.0.3}{69}{}{thm.7.0.3}{}}
\newlabel{yoneda:eq:compat.sq.1}{{7.0.4}{71}{Appendix: proof of Theorem~\ref {thm:cart.fib.chars}}{equation.7.0.4}{}}
\newlabel{yoneda:eq:defn.muhat}{{7.0.5}{71}{Appendix: proof of Theorem~\ref {thm:cart.fib.chars}}{equation.7.0.5}{}}
\newlabel{yoneda:tocindent-1}{0pt}
\newlabel{yoneda:tocindent0}{15.01021pt}
\newlabel{yoneda:tocindent1}{20.88377pt}
\newlabel{yoneda:tocindent2}{34.49158pt}
\newlabel{yoneda:tocindent3}{0pt}